\documentclass[leqno,onefignum,onetabnum]{siamltexmm}

\title{Quantification of airfoil geometry-induced aerodynamic uncertainties - comparison of approaches \thanks{This research has been conducted within the project MUNA under the framework of the German Luftfahrtforschungsprogramm funded by the Ministry of Economics (BMWi).  The authors thank   Bernhard Eisfeld,  Stefan G\"ortz and Normann Krimmelbein in Institute of Aerodynamics and Flow Technology of German Aerospace Center (DLR) for their kind help on the CFD test case. } 
} 

\author{Dishi Liu\thanks{German Aerospace Center (DLR), Institute of Aerodynamics and Flow Technology,  
              Lilienthalplatz 7, 38108 Braunschweig, Germany.
(\email{dishi.liu@gmail.com}). Questions, comments, or corrections
to this document may be directed to that email address.} 
\and Alexander Litvinenko\thanks{SRI-UQ Center, King Abdullah University of Science and Technology, 
           4700 KAUST, M.B. 4436, Thuwal,  Jeddah 23955-6900, Saudi Arabia.
(\email{alexander.litvinenko@kaust.edu.sa}).}  
\and Claudia Schillings\thanks{Seminar for Applied Mathematics ETH,
           ETH, 8092 Zurich, Switzerland.
(\email{claudia.schillings@math.ethz.ch}).} 
\and Volker Schulz\thanks{Professor, Department of Mathematics, Universit\"at Trier  
           Universit\"atsring 15, 54296 Trier, Germany.         
(\email{volker.schulz@uni-trier.de}).} }

\usepackage{graphicx}        
\usepackage{color,xcolor}
\usepackage{amssymb,amsmath,amsfonts,mathrsfs}
\usepackage{placeins}

\usepackage{epsfig}
\usepackage{epstopdf} 
 
\usepackage{float}
\usepackage{array}
\usepackage{booktabs}
\usepackage{float}

 \newcommand{\gaussfield}{\vek \psi}  
\newcommand{\vek}[1]{\mathchoice{\displaystyle\boldsymbol#1}
{\textstyle\boldsymbol#1}{\scriptstyle\boldsymbol#1}
{\scriptscriptstyle\boldsymbol#1}} 
\newcommand{\mat}[1]{\mathchoice{\displaystyle\mathbf#1}
{\textstyle\mathbf#1}{\scriptstyle\mathbf#1}
{\scriptscriptstyle\mathbf#1}}

\def\MYPsi{\textbf{\textit{S}}}

 \newcommand{\C}{\mathcal}

\begin{document}
\maketitle
\newcommand{\slugmaster}{%
\slugger{juq}{xxxx}{xx}{x}{x--x}}

\begin{abstract}
Uncertainty quantification in aerodynamic simulations calls for efficient numerical methods  to reduce computational cost, especially for the uncertainties caused by random geometry variations which involve  a large number of variables.  This paper compares five methods, including quasi-Monte Carlo quadrature, polynomial chaos with coefficients determined by sparse quadrature and gradient-enhanced version of kriging, radial basis functions and point collocation polynomial chaos, in their efficiency in estimating statistics of aerodynamic performance upon random perturbation to the airfoil geometry which is parameterized by independent Gaussian variables.   The results show that gradient-enhanced surrogate  methods achieve better accuracy than  direct integration methods with the same computational cost.
\end{abstract}



	 
\pagestyle{myheadings}
\thispagestyle{plain}
\markboth{D. Liu, A. Litvinenko, C. Schillings, V. Schulz }{Quantification of geometry-induced aerodynamic uncertainties}

\section{Introduction}
\label{sec:Intro} 
 
Aerodynamic simulations are subject to uncertainties in aircraft geometries due to various unpredictable perturbations, e.g.  manufacturing tolerances, operational wear and surface icing etc \cite{Evans07}. 
These uncertainties have the potential to dramatically lower  aerodynamic performance. Therefore enabling uncertainty quantification   is crucial for the robustness of aircraft designs.  

Uncertainty quantification (UQ) is a vivid, fast growing research field in the CFD community. 
In the context of aerodynamic design, most UQ papers focus on the treatment of uncertainties in  flow conditions, e.g.  Mach number or   angle of attack, see \cite{ResminiPhD, Bijl1,Schillings, Lucor08}, and only very little attention has been paid to uncertainties in the geometry. 

 

The quantification of geometry-induced uncertainties is computationally challenging due to large number of variables so that it  calls for efficient numerical methods. Methods sample on a regular grid are usually less favoured due to two reasons. The first is  their lack of tolerance to sample failures which are not rare in CFD computations. The second is their vulnerability to ``curse of dimensionality'', though sparse grid techniques like those employed in \cite{VSCS12,Litvin_Creta10,Resmini16} can alleviate the  curse  to some extent. Probabilistic collocation method based on Gaussian quadrature is used in \cite{Bijl} and multivariate Gaussian quadrature in \cite{Chen2011} to quantity airfoil geometric uncertainties, both with a small number of variables due to the curse. On the other side, methods based on scattered sampling schemes provide the robustness against sample failures and more freedom in choosing sample number $N$. 

Monte-Carlo (MC) and quasi-Monte Carlo (QMC) quadrature are direct integration methods using scattered sample scheme.  The MC quadrature and its variance-reduced kin (e.g.\ Latin Hypercube)  have dimension-independent error convergence rate $O(N^{-\frac{1}{2}})$,  and QMC has the worst case rate  $O(\log^d(N)N^{-1})$ with $d$ the number of variables. Due to higher degree of sample uniformity the latter can be  more accurate even with a large $d$.  
   
Surrogate-based methods are gaining more attention recently, e.g.  radial basis functions or kriging in \cite{liu2014,bompard2010,Loeven2007,giunta2004} and  polynomial chaos expansion in \cite{Kim2006,dolgov2014computation,LitvSampling13,Litvin11PAMM,Resmini16}. A non-intrusive Galerkin way of constructing surrogates is proposed in \cite{Giraldi2014}. \cite{giunta2004} shows that a kriging surrogate is better than plain  Monte Carlo and Latin Hypercube methods in estimating  mean value of a bivariate Rosenbrock function. In a recent paper \cite{Resmini16} the authors made a comparison of two classes of sample-based polynomial surrogates, generalized polynomial chaos in its sparse pseudospectral form and sparse grid stochastic collocation method, in a test case of uncertain airfoil geometry. The authors also proposed an improved version of the latter method with dimension adaptivity.

Surrogate-based methods seem to be of advantage when the gradients can be sampled at a relatively lower cost, because  direct integration methods like MC and QMC cannot effectively utilize the gradients (augmenting samples generated by Taylor's expansion are not statistically independent hence bring little benefit to the accuracy of the integration). In \cite{Kim2006} a point-collocation polynomial chaos surrogate is built to facilitate shape optimization, reinforced by gradient information acquired by local sensitivity analysis which costs at most 20\% of that for response quantities. In case of CFD model, even cheaper gradients can be obtained by using an adjoint solver \cite{brezillon2005}. By the adjoint solver all the partial gradients can be sampled at the cost of about one CFD model evaluation, i.e. each gradient costs only $1/d$ of that for response quantities.  
 
However in searching the literatures on quantifying geometry-induced uncertainties we have the following impressions: 
\begin{remunerate}
 \item In this field only one efficiency comparison of various UQ methods (not only versus Monte-Carlo method) is found in the forementioned paper \cite{Resmini16} where the authors focus in surrogate methods based on sparse grid stochastic collocation and generalized polynomial chaos.     
 \item Few papers address this problem with gradient-employing surrogate method which we consider the most efficient.
 \item No effort is made to base a UQ method comparison on  ``accurate enough'' reference statistics.  
\end{remunerate}

In this paper, we compare the efficiency of various UQ methods  in quantifying performance uncertainties  caused by geometric uncertainties of airfoil in two test cases. In the first test case the comparison includes five methods consisting of one direct integration, namely quasi-Monte Carlo (QMC) quadrature, and four surrogate-based methods: polynomial chaos (PC)  with the coefficients computed by a sparse Gauss-Hermite (SGH) quadrature (as introduced in  \cite{Barthelmann2000}), gradient-enhanced kriging (GEK), gradient-enhanced radial basis functions (GERBF) and gradient-enhanced polynomial chaos (GEPC) method (see method descriptions in Section \ref{sec:methods}).  To make the comparison more reliable, the  reference statistics (against which the error of estimated statistics is judged) is obtained by an integration of four million CFD samples and its accuracy justified sufficient by using multipartition method. In the second test case three methods are compared, namely QMC, plain kriging and GEK. The number of CFD samples utilized by the UQ methods is kept small ($ \leq 500$) in the both test cases to make the comparisons more meaningful to industrial applications. 
  
For the geometric uncertainties in the test cases, we focus on smooth geometry deviations from the original shape as those typically  arise during the manufacturing process.  Other  uncertainties or errors, such as those due to physical, mathematical and numerical modeling approximations, are not considered in this work.    

The plan for the rest of the paper is as follows: Section \ref{sec:methods} briefly introduces the methods in the comparison.  Section \ref{sec:test case 1} sets up a test case based on an Euler model of an airfoil with uncertain geometry, and gives results of the comparison as well as a discussion. Section \ref{sec:test case 2} does the same for a Navier-Stoke test case.   Section \ref{sec:conclusion} gives a summary of the paper.

\section{Methods}\label{sec:methods} 

In this section the methods in the comparison are introduced. More weights are given to  two novel implementations, the  gradient-enhanced radial basis functions  and gradient-enhanced  polynomial chaos    method. Only brief introductions or references are given to the others.  
 
For the  Quasi-Monte Carlo (QMC) quadrature \cite{caf1998}, Sobol low discrepancy sequence generated by the algorithm given in \cite{joe2008} is used since the sequence is shown comparable or slightly better than other sequences  in a numerical comparison \cite{radovic1996}. 

The  gradient-enhanced kriging (GEK) and plain kriging \cite{chung2002} are implemented  by using  \textit{Surrogate-Modeling for Aero-Data Toolbox (SMART)} \cite{han2012} developed at DLR, opting for  ordinary kriging and a correlation model of cubic spline type which is considered  the most efficient in similar situations in \cite{liu2012b}.  The internal parameters of the correlation model are fine-tuned to fit the sampled data by a maximum likelihood estimation \cite{zimmermann2010}.   

The gradient-enhanced radial basis function  method is introduced in  Section \ref{subsec: GERBF},   polynomial chaos  method based on a sparse Gauss-Hermite quadrature (PC-SGH) and gradient-enhanced  polynomial chaos method (GEPC) are introduced in Section \ref{subsec: GEPC}.
  
To quantify uncertainty via surrogate methods, one first establishes the surrogate based on a small number of QMC samples of the CFD  model, and integrates for the target statistics and probability density function (pdf) by a large number  (we take $1 \times 10^6$) of QMC samples on the surrogate, only except for the polynomial chaos methods where the mean and variance can be directly obtained.

 \subsection{Gradient-enhanced radial basis functions (GERBF)} \label{subsec: GERBF} 
  Denoting the CFD model as $f$, a radial basis function (RBF) \cite{buhmann2000} approximate takes the form  
\[
  \hat f(\vek \xi) = \sum_{i=1}^N  w_i \, \phi_i(\| \vek \xi- \vek \xi^{\langle i \rangle}\|),
\]
where $\phi_i$  are radial basis functions, $\| \cdot \|$ denotes the Euclidean norm, and $\vek \xi^{\langle i \rangle}$ is a sample point where  $\phi_i$ is radial about.  The coefficients $w_i$  are determined  by  fitting    $\hat f(\vek \xi)$ to $N$ samples.

Denoting the Euclidean distance from the center as $ r$, popular types of $\phi( r)$ include  $\sqrt{ r^2 + a^2}$ (multiquadric),  $1/\sqrt{ r^2 + a^2}$ (inverse multiquadric),  $\exp (-a^2  r^2)$ (Gaussian) and $ r^2 \ln (a r)$ (thin plate spline), in which $a$ is a parameter to be fine-tuned for a particular set of samples. Gradient-assisted RBF were proposed in \cite{gian2006,ong2008} where the gradients are exploited in the same way as gradient-enhanced kriging (GEK).  This method is numerically equivalent to a GEK that uses   RBF as the correlation function.

A different gradient-employing RBF method  is implemented in this work, which only requires RBF to be first-order differentiable, in contrast to the second-order differentiability required by the gradient-assisted RBF. We call it gradient-enhanced RBF (GERBF) in this paper. To accommodate the gradient information,  this method introduces additional  RBF that are centred at \textit{non-sampled} points, i.e. an GERBF  approximation is
\[
  \hat f(\vek \xi) = \sum_{i=1}^{N(1+d)}  w_i \, \phi_i(\| \vek \xi- \vek \xi^{\langle i \rangle}\|), \;\;\;\;\; \mbox{with $d$ the number of variables  }\; \;\;  
\]
The $\vek \xi^{\langle i \rangle}$ with $i \leq N$ are sampled   points, while those with $i>N$ are non-sampled  points which can be chosen arbitrarily as long as none of them duplicates (or too close to) the sampled ones. The coefficients $\vek w = \{  w_1, \cdots, w_{N(1+d)}\}^T$  are determined by   fitting   $\hat f(\vek \xi)$ to $N$ samples and $Nd$ gradients,   i.e. solving the following system,  
 \[
   \mat A \vek w = \vek f .
 \]
Abbreviating $\phi_i(\| \vek \xi  - \vek \xi^{\langle i \rangle}\|)$ as $\phi_i(\vek \xi )$, $\mat A$ is a $N(1+d)$-sized square matrix whose first $N$ rows contains   values of $\phi_i(\vek \xi)$ evaluated at the $N$ sampled points and the rest rows   values of $\partial \phi_i/ \partial \xi_j$ for $1 \le j \le N $ evaluated at the $N$ sampled points.  $\vek f$ is a $N(1+d)$-length vector containing the sampled values of $f$ and its $d$ partial derivatives.

Since $\mat A$ is often ill-conditioned the system is solved through a truncated singular value decomposition \cite[Chapter 15]{Press2007} in which singular values smaller than a threshold are discarded\footnote{Let $\Lambda_i$ be the singular values, those $|\Lambda_{i}|/ \max_i (| \Lambda_{i}|)\leq N(1+d)\times 10^{-13} $ are discarded.}.
This gives a solution with the smallest $L_2$ norm which is less spoiled by the ill-conditionedness \cite[Chapter 15]{Press2007}.  For  large systems a more efficient and stable solution by pivoted Cholesky decomposition is proposed in \cite{liu2015}.

Inverse multiquadric RBF is used in this comparison. The internal parameter $a$ is fine-tuned by  a \textit{leave-one-out} cross-validation as in \cite{bompard2010}.  
 
\subsection{Polynomial chaos  methods}  \label{subsec: GEPC} 
 
According to Wiener \cite{wiener1938},  a surrogate could also be made in the form of a truncated {\em  polynomial chaos expansion} (PCE) :
\begin{align}
 \hat  f(\vek \xi) = \sum_{i=0}^{K} c_i H_i (\vek \xi)  \label{PCE}
\end{align}
where $H_i$ is a multivariate Hermite polynomial  with Gaussian variables (see introduction in e.g. \cite{matthies2007}). It is called ``truncated'' because the $K$ is taken finite for practical computation. $K$ is related to the maximal order of polynomial ($p$) and number of variables ($d$) by $K= (p+d)!/(p!d!)$.  

The coefficients $c_i$ in Equation (\ref{PCE}) could be determined by various ways\cite{le2010spectral,xiu2010numerical}, in this work the following two of them are implemented.
\subsubsection{Projection based on sparse Gauss-Hermite quadrature (SGH)}  
Due to the mutual orthogonality of $H_i$'s with respect to their inner product,  the coefficients $c_i$ can be  expressed by \[
   c_i = \frac{\int f(\vek \xi) H_i(\vek \xi) \rho(\vek \xi) \mbox{d} \vek \xi}{\int H_i^2(\vek \xi) \rho(\vek \xi) \mbox{d} \vek \xi},   
  \]  
 where $\rho(\vek \xi)$ is  Gaussian probability density. The numerator is approximated by using sparse Gauss-Hermite (SGH) quadrature, the denominator is computed analytically\cite{matthies2007}. 
%
\subsubsection{Gradient-enhanced  point-collocation  polynomial chaos  method (GEPC)}  
Point-collocation  polynomial chaos  method \cite{hosder2006} determines   $c_{i}$ by fitting equation (\ref{PCE}) to $N$ collocation points. The work \cite{hosder2007} suggests to use $N=2K$ for the best performance   which leads to an over-determined system   to be solved by a least square approach.  A gradient-employing version of this method is implemented in this work, in the name gradient-enhanced  polynomial chaos (GEPC) method, where accordingly one chooses $N(1+d)=2K$. The unknown coefficient $\vek c = \{c_0, c_1, \cdots, c_K \}^T$ is determined by solving the following system  
 \[
   \mat \Psi  \vek c = \vek f 
 \]
in which  

\begin{math}
\mat \Psi  = 
\left [  
\arraycolsep=1.4pt\def\arraystretch{0.5}
\begin{array}{c} 
            \mat \Psi_0   \\        
             \vdots     \\        
               \mat \Psi_d          
            \end{array}
\right ], \;\;\;  
\mat \Psi_0  = 
\left [  
\arraycolsep=1.4pt\def\arraystretch{0.5}
\begin{array}{lllll} 
            H_0 (\vek \xi^{\langle 1 \rangle}) &  \;\; &  \cdots & \;\;  & H_{K} (\vek \xi^{\langle 1 \rangle})    \\            
           \hspace{4mm} \vdots   &  \;\;    &  \ddots & \;\;  &  \hspace{4mm}\vdots    \\   \vspace{0.2cm}
            H_0 (\vek \xi^{\langle N \rangle}) &  \;\;  &  \cdots & \;\;  & H_{K} (\vek \xi^{\langle N \rangle})                       
            \end{array}
\right ], \;\;
 \end{math} 
  \vspace{0.2cm} 
  
 \begin{math}
\underset{d\geq j>0}{\mat \Psi_j}  = 
  \left [ 
 \arraycolsep=1.4pt\def\arraystretch{0.5}  
   \begin{array}{lllll} 
            H^{(j)}_0 (\vek \xi^{\langle 1 \rangle}) &  \;\;  &  \cdots & \;\;  & H^{(j)}_{K} (\vek \xi^{\langle 1 \rangle})    \\      
            \hspace{4mm}\vdots   &  \;\;  &  \ddots & \;\;  &\hspace{4mm} \vdots    \\     
            H^{(j)}_0 (\vek \xi^{\langle N \rangle}) &  \;\; &  \cdots & \;\;  & H^{(j)}_{K} (\vek \xi^{\langle N \rangle})    \\                                
            \end{array}
\right ]
 \end{math} \vspace{0.3cm}\\  
 with $H^{(j)}_i = \partial H_i / \partial \xi_j$. Thus, $\mat \Psi$ consists of values of the  $H_i$ and   $\partial H_i / \partial \xi_j$ evaluated at $N$ sample points, sized $N(1+d) \times K$. The vector $\vek f$ is composed as in the GERBF method.     
 

 
Once a surrogate based on polynomial chaos is established,  the estimates of the mean and  the variance  of  $f(\vek \xi) $  can be directly obtained  from $\vek c$:
\begin{align}
\mu =  c_{0}\; , \;\;\;\;\;\;  \sigma^2 \approx  \sum^{K}_{i=1} \left [ c^2_i    \int H^2_i(\vek \xi) \rho(\vek \xi)\mbox{d} \vek \xi \right], \label{variance}
\end{align}
where the integration can be made analytically. The exceedance probabilities and pdf are integrated by a large number (we take $1 \times 10^6$) of QMC samples on the surrogate model.

\section{An Euler test case of uncertain airfoil geometry}
\label{sec:test case 1}

\begin{figure}[htbp]
  \centering
  \begin{minipage}[b]{0.49\textwidth}
\centering
    \includegraphics[width=\textwidth]{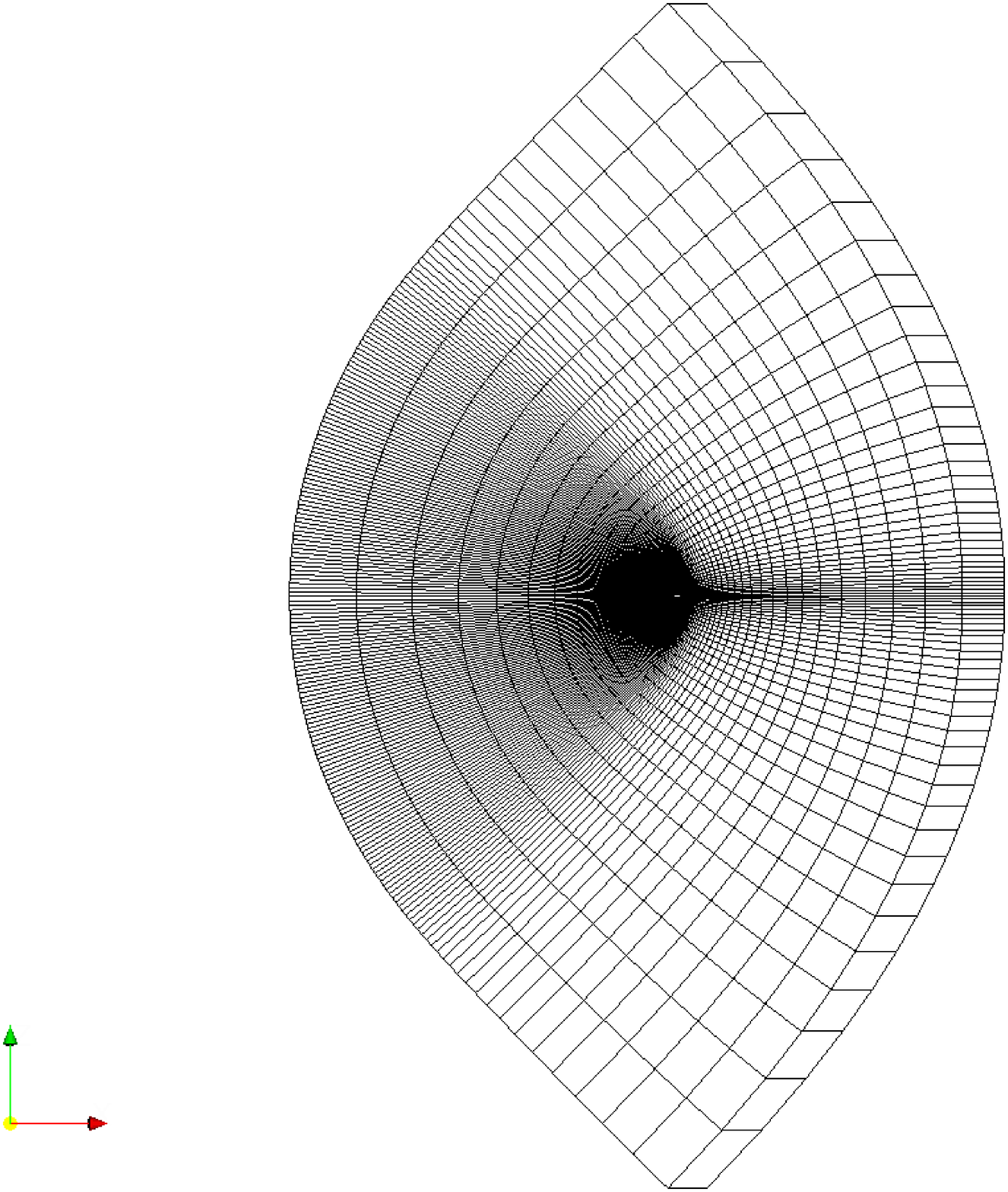}  
  \end{minipage}
  \begin{minipage}[b]{0.49\textwidth}
\centering
    \includegraphics[width=1.0\textwidth]{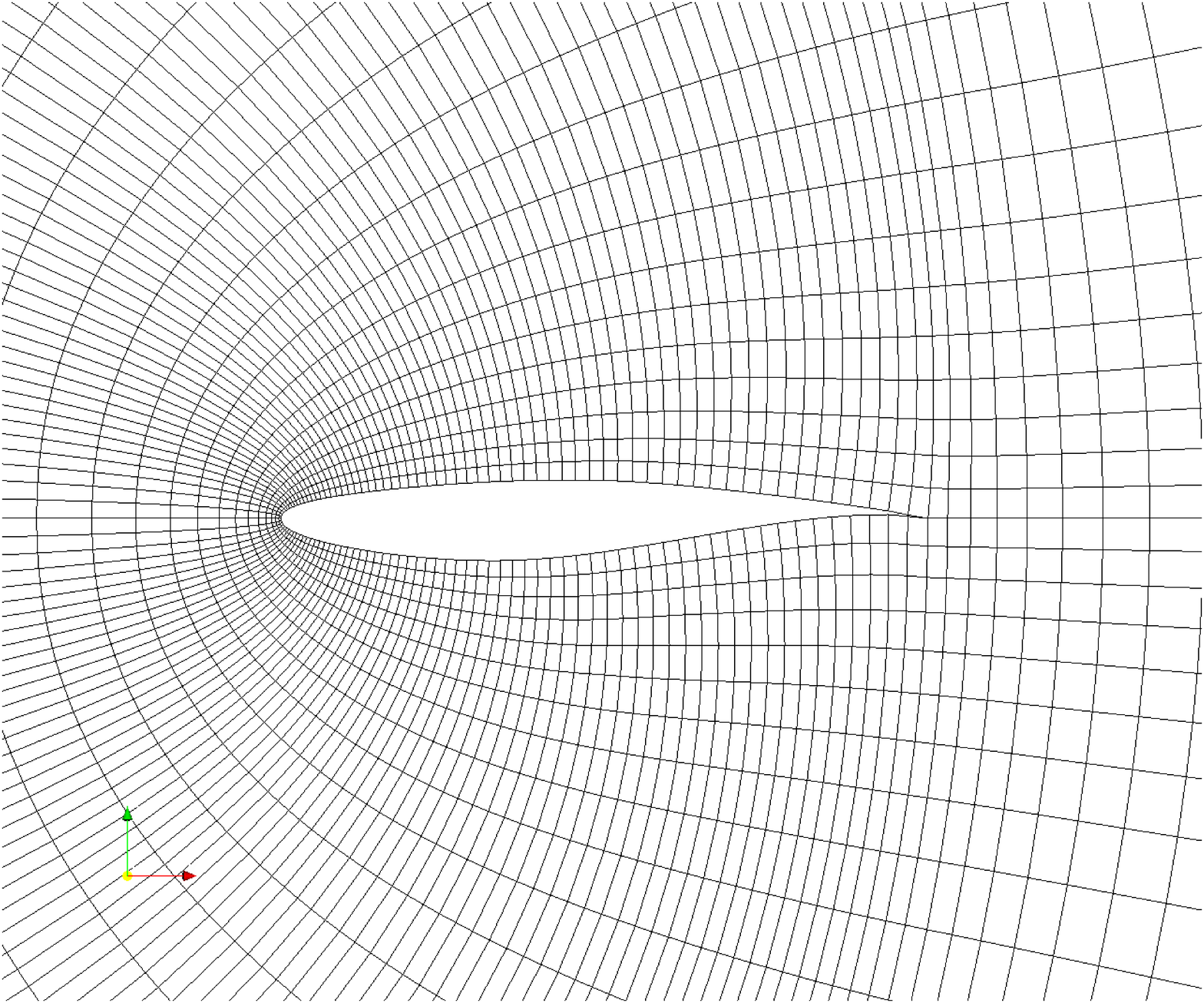}  
  \end{minipage}
\caption{\label{gridFlowerRAE} Grid for the RAE2822 airfoil in the Euler test case: the total grid (left) and zoom around the airfoil (right).}
\end{figure}

The first test case   is a CFD model of inviscid flow around a  RAE2822 airfoil at Mach number  0.73, angle of attack (AoA) 2.0 degree   and temperature 273.15 degree.
We use flow solver TAU \cite{TAU2}, opting for central flux discretization, scalar dissipation, backward Euler solver and  ``4w''  multigrid cycle.   The domain is discretized by a $193\times33$ grid in which the airfoil has 128 surface nodes, as shown in Figure~\ref{gridFlowerRAE}.

\begin{figure}[htbp]
    \includegraphics[width=0.7\textwidth]{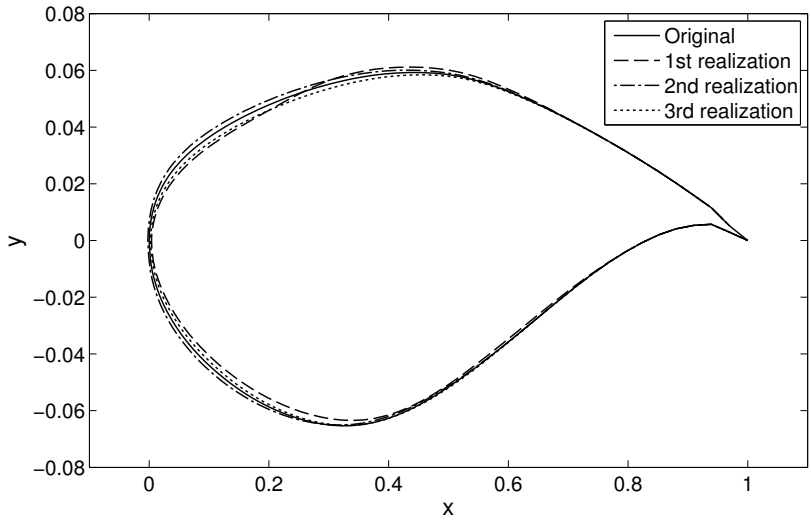}  
\caption{\label{resulpertshapes} Three realizations of  perturbed geometry}
\end{figure}

The source of uncertainty is a random perturbation to the original airfoil geometry, which is modeled by a random field. Representing the original geometry (see Figure~\ref{resulpertshapes}) by a matrix $\mat G \in \mathbb R^{128\times 2}$ which  contains  the $x,y$ coordinates of the surface nodes,  the random field of perturbation can be written as $\vek R(\vek g, \omega) \in \mathbb R^{128}$ with $\vek g   =\{x, y\}  \in \mat G$ and $\omega \in \Omega$,  $\Omega$ is the sample space of the random perturbation. 
 
$\vek R(\vek g, \omega) $ is defined as a transformation of a Gaussian random field $\vek \psi(\vek g, \omega)$ which can be easily parameterized by independent random variables by using Karhunen-Lo\`{e}ve expansions.  

The Gaussian field $\vek \psi(\vek g, \omega)$ is defined by a zero mean and  a Gaussian-type correlation function:
 \begin{equation}\label{cov}
 \mbox{cov}\left(\vek g_i, \vek g_j\right)=\sigma(\vek g_i)\ \sigma(\vek g_j) \exp \left(-\frac{\|\vek g_i- \vek g_j \|^2}{\ell^2}\right)
\end{equation}
 with  correlation length $\ell=0.005$, the standard deviation  $\sigma$ is assumed to be: 
   \[  \sigma(\vek g)= \left \{ \begin{array}{ll}                              
                                (0.8-x)^{0.75} & \mbox{  if  } x \leq 0.8  \\ 
                              0  & \mbox{  otherwise } 
                   \end{array}
                   \right.
                   \] 
The zero perturbation to the trailing end is to guarantee the convergence of the CFD solver.
 
To ensure  boundedness of the perturbation,  the   random perturbation field  $\vek R$ is made to have a sine-shaped probability distribution by such  a transformation of the Gaussian field $\vek \psi$ :
\begin{equation}
\vek R(\vek g, \omega)=2 s(x)\cdot \arccos \left(1-2\Phi \left(\vek \psi\left(\vek g,\omega\right)\right)\right)/ \pi-1, \label{transGRF}
  \end{equation}
where $\Phi$ is the Gaussian cumulative distribution function and  $s(x)=(0.8-x)^{0.75}\cdot\sqrt{0.00002}$.

Perturbed geometries are generated by imposing the perturbation $\vek R$ to the original geometry, in the direction normal to the airfoil surface:   
 \[
  \widetilde {\mat G}  (\omega) = \mat G +  \vek R(\vek g, \omega)\cdot \vek n(\vek g ),
 \] 
where $\vek n(\vek g)$ is the normal vector at $\vek g$. 

However  the above representation has 128 correlated variables  while one prefers a representation with independent variables in a fewer number.   A model reduction can be made by an approximation of $\vek \psi$ through a Karhunen-Lo\`{e}ve expansion (KLE)  \cite {adler2007}, i.e.
\begin{equation}
 \gaussfield (\vek g,\omega) \approx \sum\limits_{i=1}^k \sqrt{\lambda_i}\; \vek V_i(\vek g)\; \xi_i(\omega) = :   \widehat{ \gaussfield}(\vek g,\omega) \label{KLE}
 \end{equation}
 where $\lambda_i$ and $\vek V_i$ denote the eigenvalues and eigenvectors of  the covariance matrix generated by Equation (\ref{cov}) and $\xi_i$'s are independent standard Gaussian random variables.  For our problem the number of variables is reduced to $k=9$, this truncated KLE approximation retains $99.98\%$ of  variance of the random field. The KLE is optimal in all linear form representations with the same number of terms since it minimizes the loss in the variance  caused by the model reduction\cite{adler2007}.  In numerical implementation   $ \widehat{  \gaussfield}(\vek g,\omega)$ is used instead of $\gaussfield (\vek g,\omega)$ in Equation (\ref{transGRF}), so that the random field of perturbation $\vek R$ is  parameterized by 9 independent Gaussian variables.  In \cite{khoromskij2009application,liu2015} efficient algorithms for eigen-decomposition of large covariance matrices are given, by which one has the choice to compute only the  $k$ largest eigen modes. 
 
\begin{figure}[htbp!]
\centering
\includegraphics[width=0.48\textwidth]{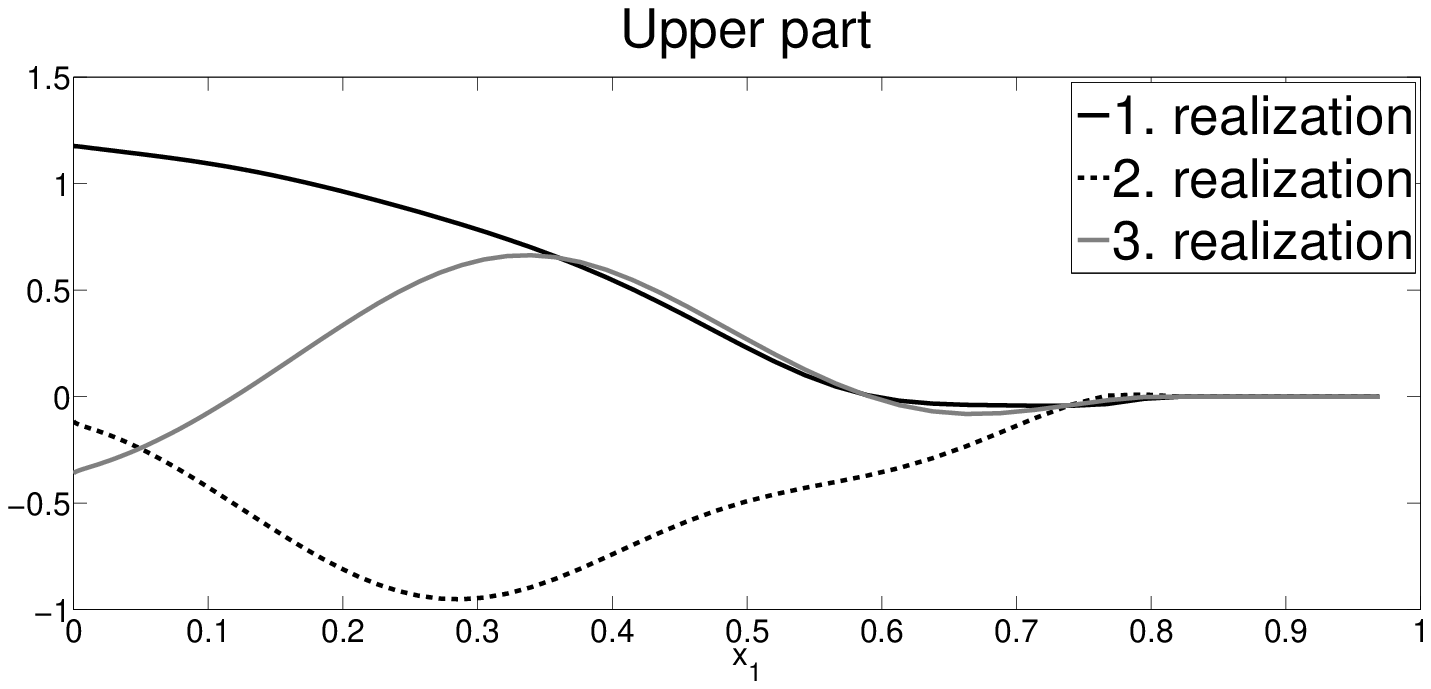}
\includegraphics[width=0.48\textwidth]{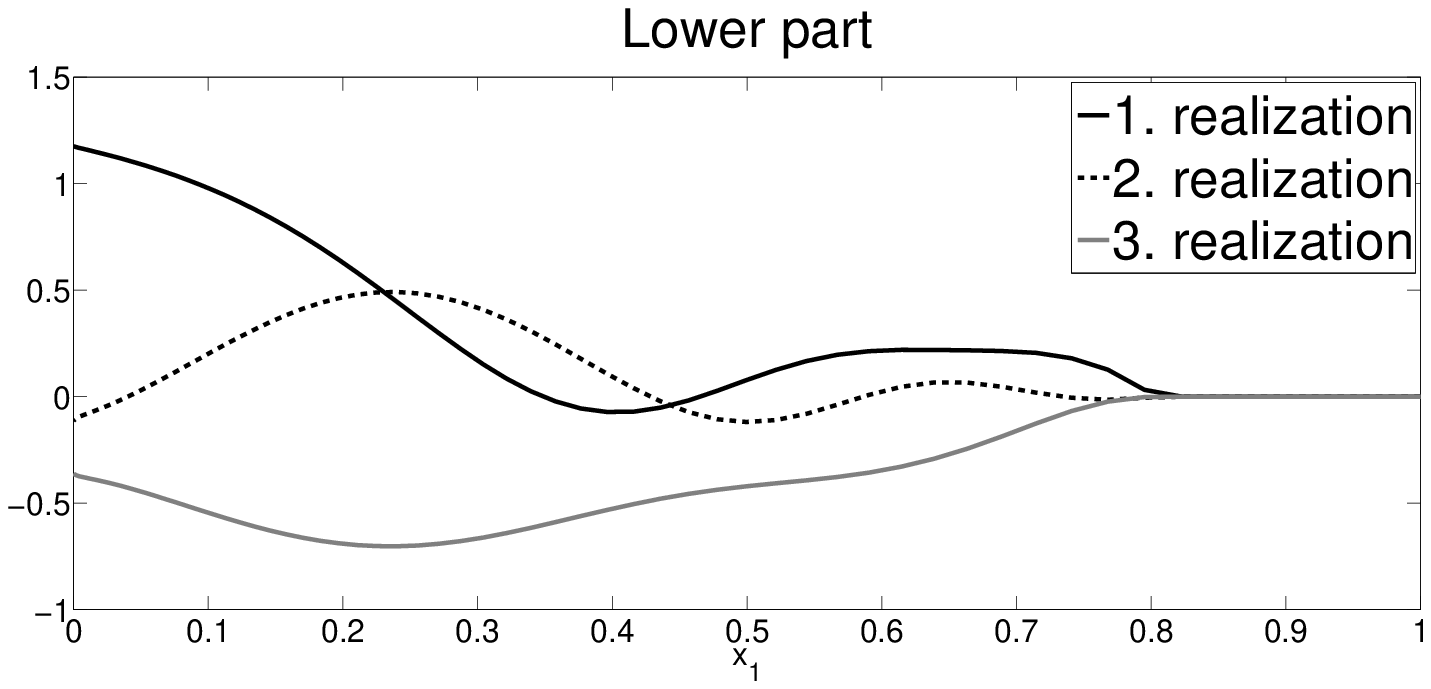} 
\caption{\label{figreal} Realizations of the Gaussian random field $\gaussfield$: upper  (left) and lower (right) part of airfoil.}
\end{figure}

\begin{figure}[htbp]
\centering
\includegraphics[width=0.48\textwidth]{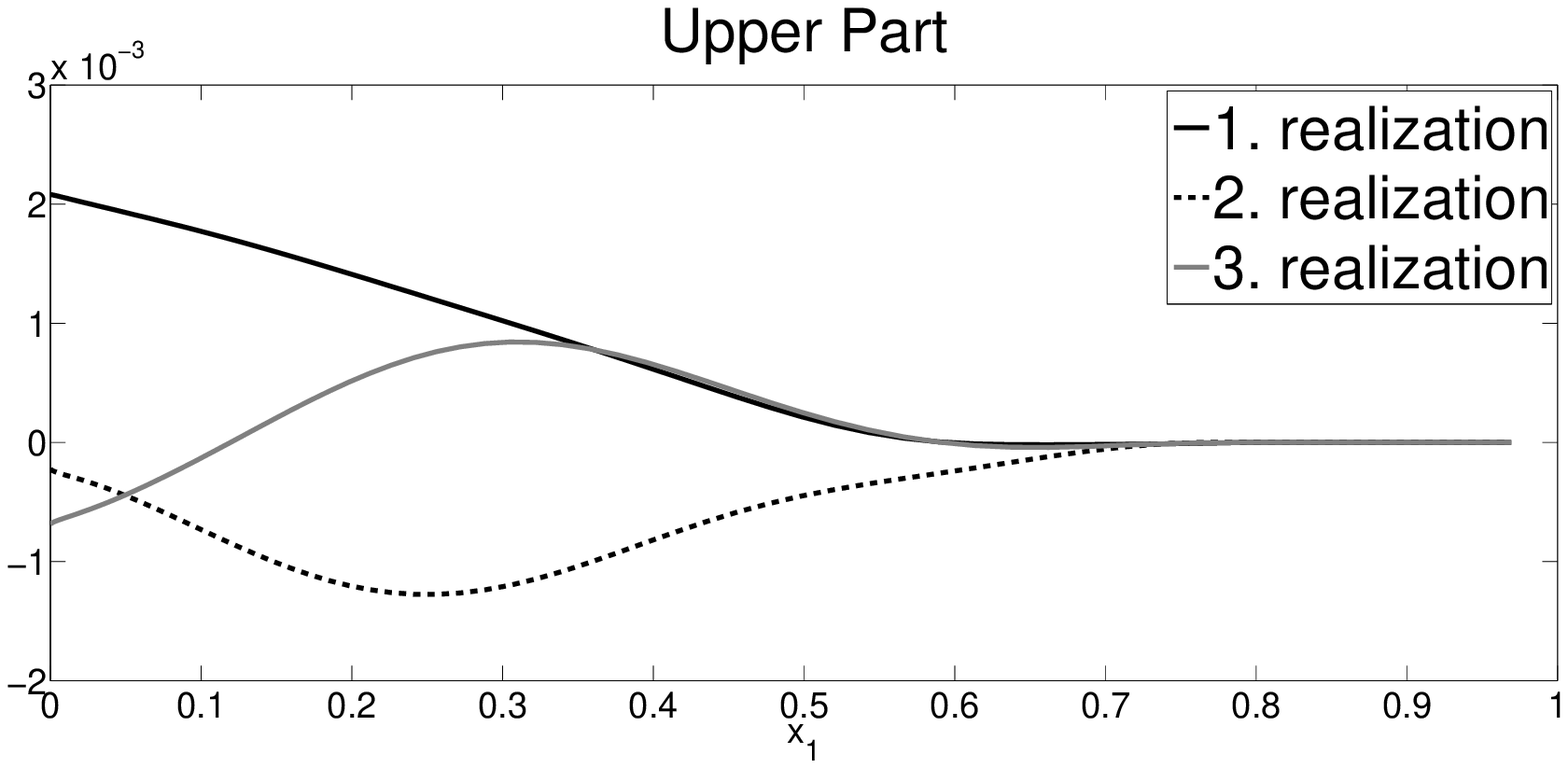}
\includegraphics[width=0.48\textwidth]{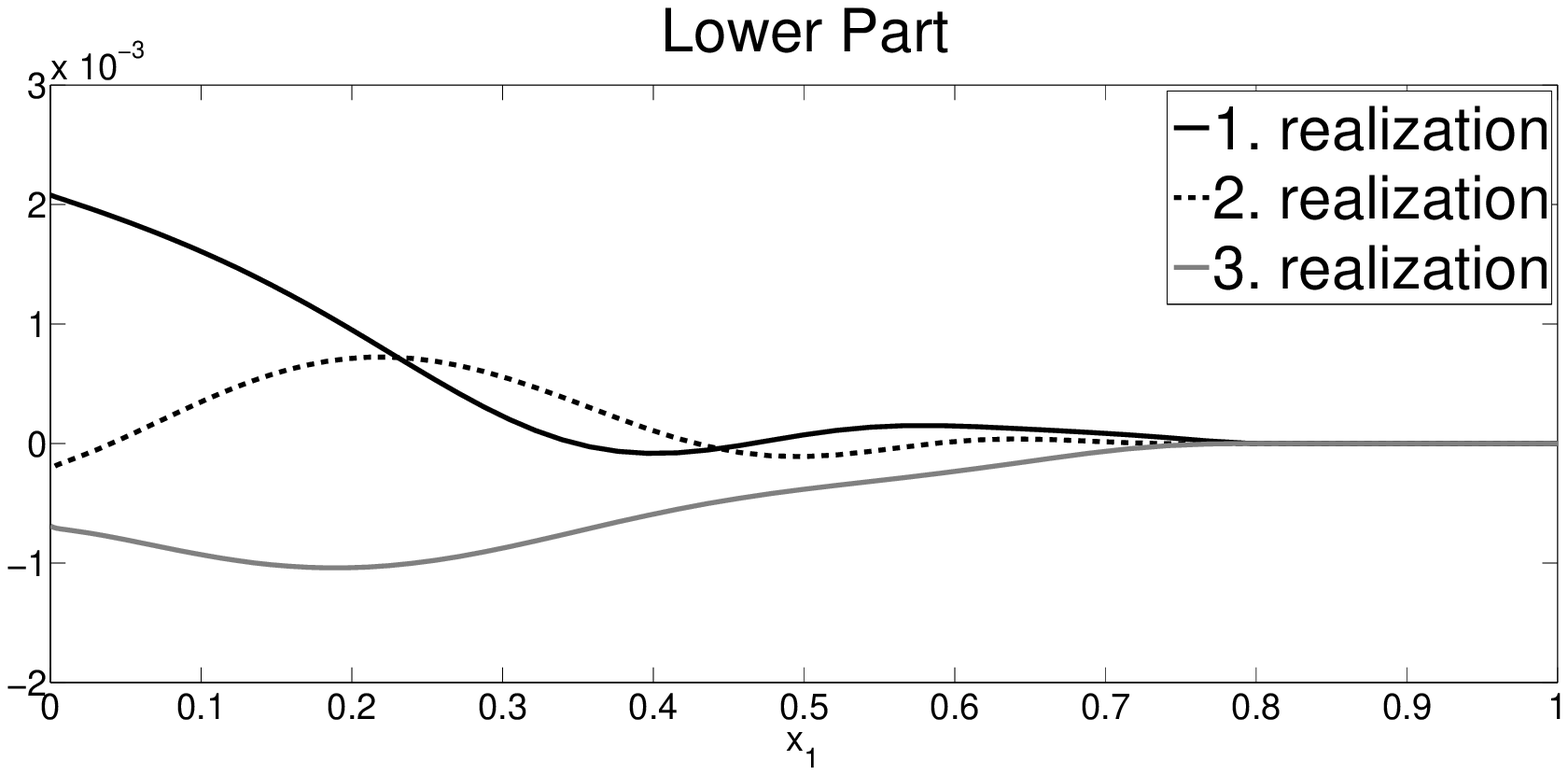} 
\caption{\label{transreal} Realizations of the transformed random field $\vek R$:  upper (left) and lower   (right) part of airfoil.}
\end{figure}

\begin{figure}[htbp]
    \includegraphics[width=0.6\textwidth]{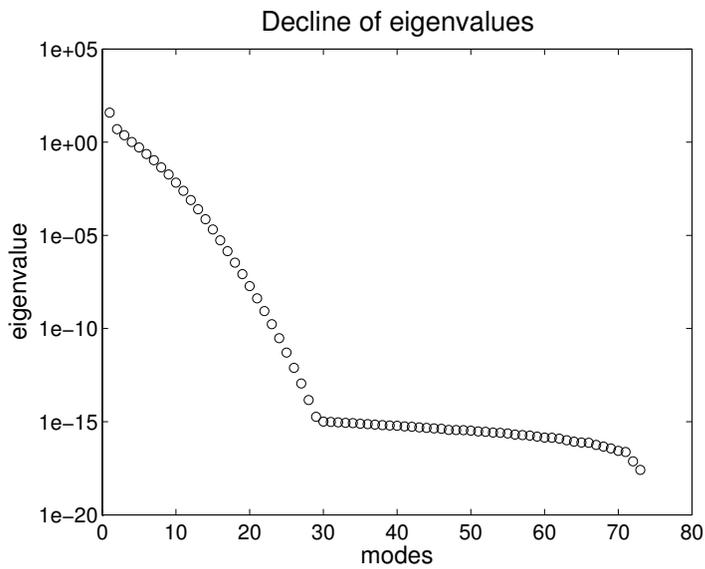}  
\caption{\label{KLeiv} Decline of eigenvalues}
\end{figure}

 
Figure~\ref{figreal} and \ref{transreal} display three  realizations of the  $\gaussfield$ and $\vek R$ respectively.  Examples of perturbed  geometries  $\widetilde {\mat G}$  are shown in Figure~\ref{resulpertshapes}. Figure~\ref{KLeiv} shows the fast decline of the eigenvalues $\lambda_i$.  

The test case is a scaled-down model with an inviscid flow assumption and relatively coarse mesh. This simplifying configuration is to facilitate a large number of evaluations to obtain reliable reference statistics of performance on the basis of which the accuracy of estimated statistics are compared. Since this accuracy is generally  independent of the model bias caused by the particular configuration, the result of this comparison holds meaningful and representative in  judging the efficiency rank of UQ methods in other stricter configurations.

 
\subsection{Settings of the comparison}\label{sec:setting}
  With the above parameterization of geometric uncertainties, the five UQ methods are applied to the test case and compare their efficiency in estimating statistics of two aerodynamic system response quantities (SRQ), namely the coefficients of lift and drag ($C_L$ and $C_D$), as well as their probability distribution functions (pdf).
\subsubsection{Target statistics and their reference values}\label{sec:reference_values}
The following statistics of  $C_L$ and $C_D$ are estimated:
\begin{itemize}
 \item   means   $\mu_L$ and $\mu_D$
 \item    standard deviations $ \sigma_{L}$ and $\sigma_D$ 
 \item    exceedance probabilities   
$P_{L, \kappa} = \mbox{Pro} \{ C_L \leq \mu_L - \kappa \cdot  \sigma_{L} \}$ \quad \text{and} $\quad  P_{D,\kappa} = \mbox{Pro} \{ C_D \geq \mu_D + \kappa \cdot  \sigma_D \}$ with $\kappa=2,3$
\end{itemize} 
 
The reference values of these statistics are obtained from a relatively large number  ($N=4\times 10^6$) of  quasi-Monte Carlo (QMC) samples of the CFD model. To justify the validity of using these reference values in the efficiency comparison, their accuracy are estimated using Snyder's multipartition method \cite{snyder2000} since the theoretical error bound of QMC integration is not a practical accuracy indicator. Suppose $\MYPsi$ denotes any of the above statistics, by this method one makes an equal-size $m$-partition of the $4\times 10^6$ samples,  and obtained
$m$ estimates $\{ \widetilde \MYPsi_i\}^m_{i=1}$ by integrating on each partition only, then compute the sample standard deviation $\varsigma_{m}$ of  $\widetilde \MYPsi_i$. $\varsigma_m$ is an estimate of the standard deviation of QMC estimate of  $\MYPsi$ with sample number $N/m$. To extrapolate $\varsigma_m$ to $\varsigma_1$ one computes $\varsigma_m$ for four values of $m$, i.e. $m=\{128, 64, 32,16 \}$, and fit a line across the  four $\left (\log(N/m), \log(\varsigma_m) \right )$ points using weighted linear least square method, the $m$ values are used as weights to account for the increasing variability for smaller values of $m$. With the fitted slope $\alpha$ and intercept $\beta$ the $\varsigma_1$ is extrapolated as  $ \varsigma_1 = \exp(\alpha \log(N) + \beta ) $
which is an estimate of the standard deviation of the QMC estimate of  $\MYPsi$ with $N$ samples.

$\varsigma_1$ values for each reference statistics are listed in Table \ref{T2}.  
 In the efficiency comparison of the UQ methods (results given in section \ref{sec:results}) the smallest  measured  errors (measured  against  these reference values) in   mean and standard deviation (stdv) are at least by 10 times larger than $3\times\varsigma_1$, which means, by taking the assumption that the reference values are Gaussian distributed around the \textit{true} values of the statistics,  the measured   errors have a $99.73\%$ confidence interval  of at widest $\pm 10\%$.  For exceedance probabilities this confidence interval is  also  valid except for a few measured errors, as shown in  Figure~\ref{fig:error_CL} and \ref{fig:error_CD} where the values of the corresponding  $3\times\varsigma_1$ are depicted by thick dash line. 
\begin{table}[htbp]    
\footnotesize
\centering
\begin{tabular}{  c|c|c|c}  
\hline  \hline  
    $\; \varsigma_1 (\mu_L) \;$   &    $ \;\varsigma_1 (\sigma_L) \;$      &    $ \;\varsigma_1 (P_{L,2}) \;$  &$\; \varsigma_1 (P_{L,3}) \; $\\  \hline  
       9.1e-9         &       5.0e-8                     &   1.8e-5                     &       9.3e-6   \\  
 \hline   \hline  
\end{tabular}    \vspace{0.4cm}\\
 \begin{tabular}{  c|c|c|c}  
\hline  \hline  
    $\; \varsigma_1 (\mu_D) \;$   &    $ \;\varsigma_1 (\sigma_D) \;$      &    $ \;\varsigma_1 (P_{D,2}) \;$  &$\; \varsigma_1 (P_{D,3}) \; $\\  \hline  
     3.0e-9         &       1.4e-8                     &   1.2e-5                     &       6.8e-6   \\  
 \hline   \hline  
\end{tabular}  
\caption{Estimated stdv  of  reference statistics for $ C_L$ (top) and $C_D$ (bottom)}
\label{T2}
 \end{table}
 

\subsubsection{Measure of cost}
 
The computational cost is  measured in term of  ``compensated evaluation number'' $N_c$.  For the three gradient-employing methods,  $N_c=3N$   with $N$ denoting the  number of CFD model evaluations,  since with an adjoint solver the computational cost of the gradients of each system response quantity (SRQ)   equals  approximately to the cost of one CFD  evaluation and we need  the gradients for  two SRQs ($C_L$ and $C_D$).  For QMC and PC-SGH methods,   $N_c=N$. 

\subsubsection{Range of sample numbers }
The number of CFD evaluations is kept smaller than 500 in this comparison, as shown in the Table \ref{T2_a}. Most of the methods share the same   range except for PC-SGH in which only quadrature level 2 and 3 are used since the sample number of level 4 is well beyond 500.
 
 \begin{table}[htbp]   
 \footnotesize
 \centering
 \begin{tabular}{ c | c | c }
 \hline \hline
    Methods &   $N$ & $N_c$ \\ \hline \hline 
   QMC &    30 $\sim$ 420 & 30 $\sim$ 420  \\ \hline  
   PC-SGH (order 2, 3)&  19 and 181 & 19 and 181  \\ \hline  
    GEK &    10 $\sim$ 140 & 30 $\sim$ 420  \\ \hline  
    GERBF &    10 $\sim$ 140 & 30 $\sim$ 420  \\ \hline  
    GEPC (order 2,3 and 4)   \hspace{0.2cm}& \hspace{0.2cm}   11, 44 and 143  \hspace{0.2cm}& \hspace{0.2cm} 33, 132 and 429  \hspace{0.2cm} \\    
\hline \hline  
\end{tabular}
\caption{ Sample number  in the comparison} \label{T2_a}
\end{table}

\subsection{Numerical results  and discussion }\label{sec:results}
 
\begin{figure}[htbp!]
  \centering
 \includegraphics[width=0.43\textwidth]{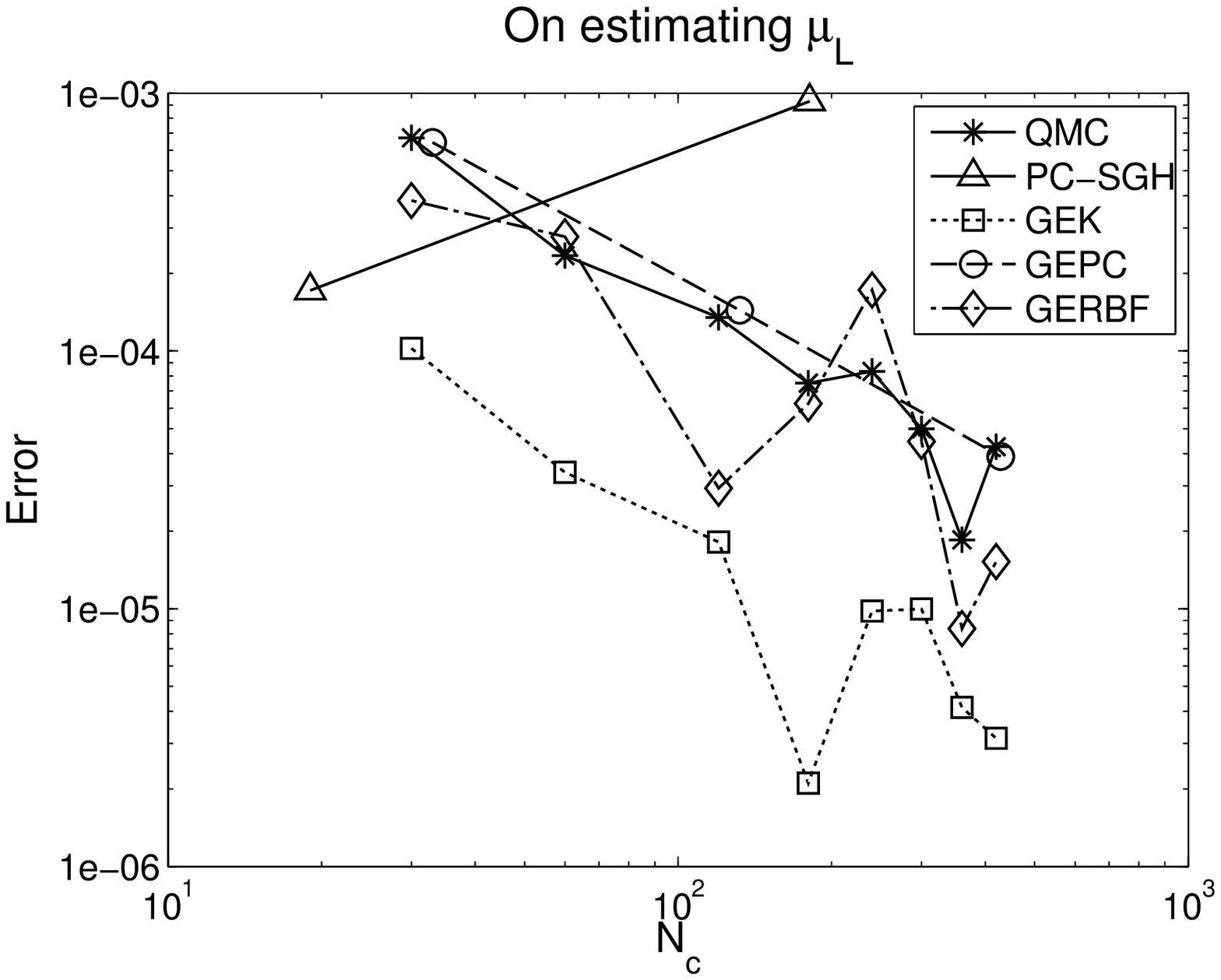}
 \includegraphics[width=0.43\textwidth]{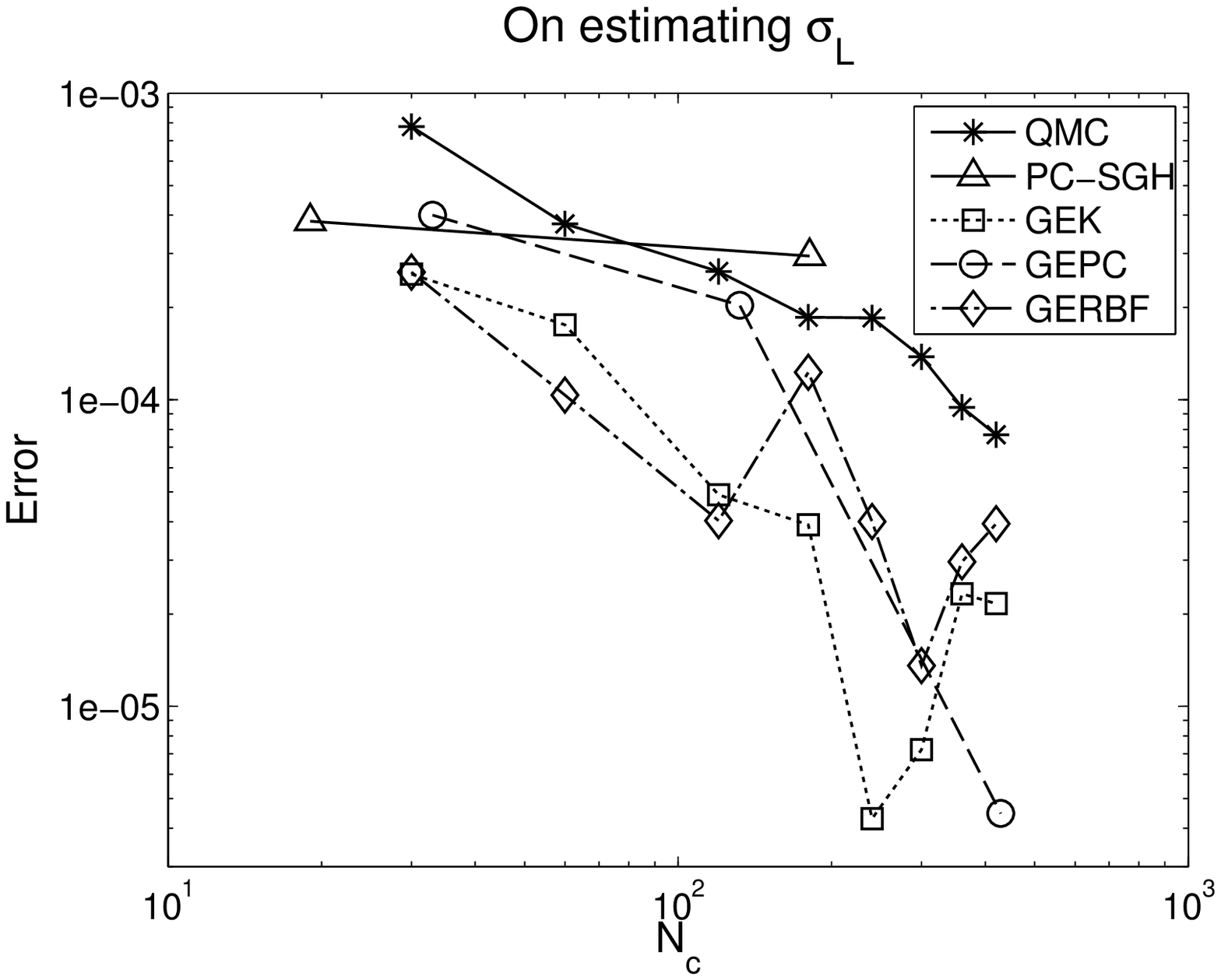}  
 \includegraphics[width=0.43\textwidth]{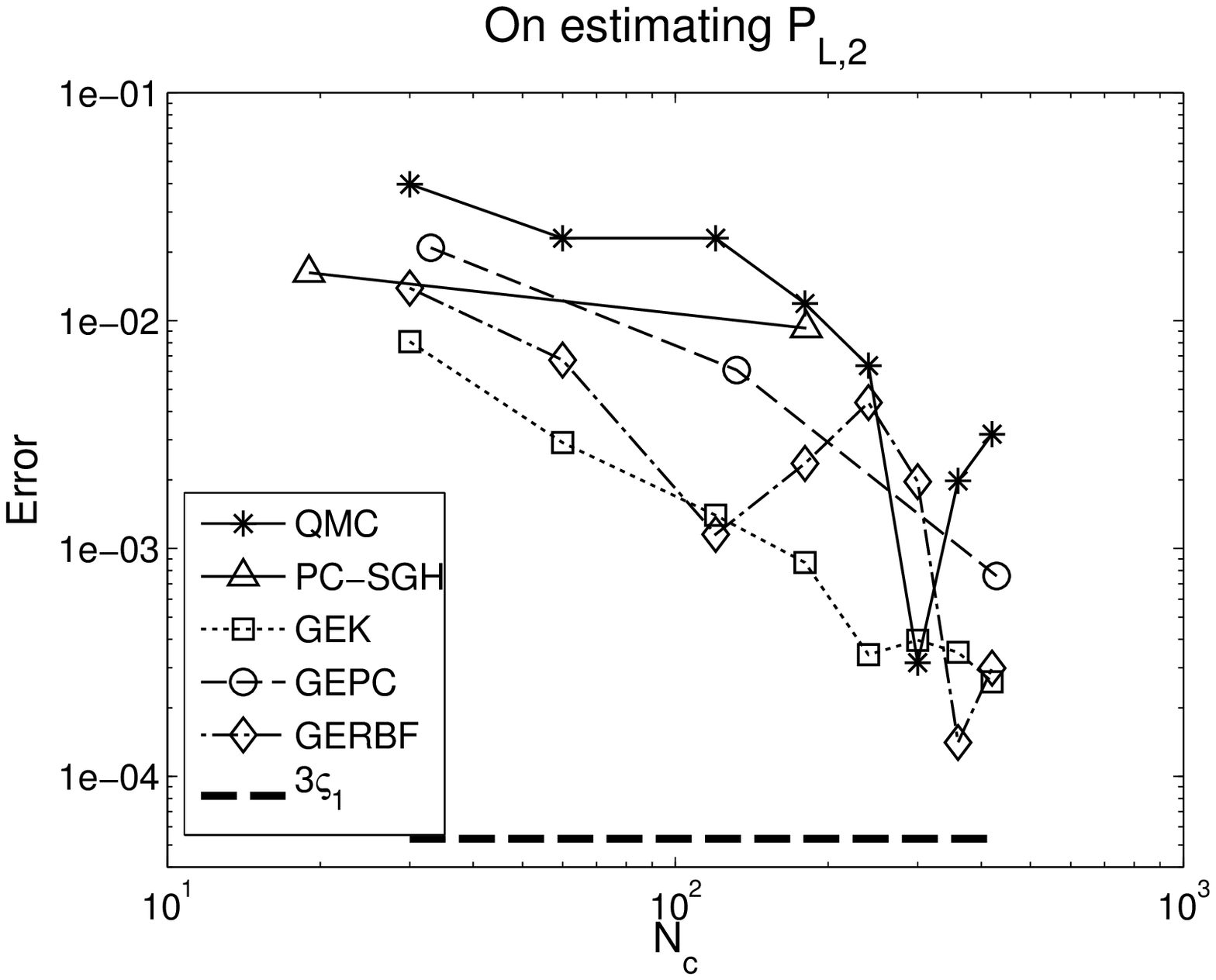}
 \includegraphics[width=0.43\textwidth]{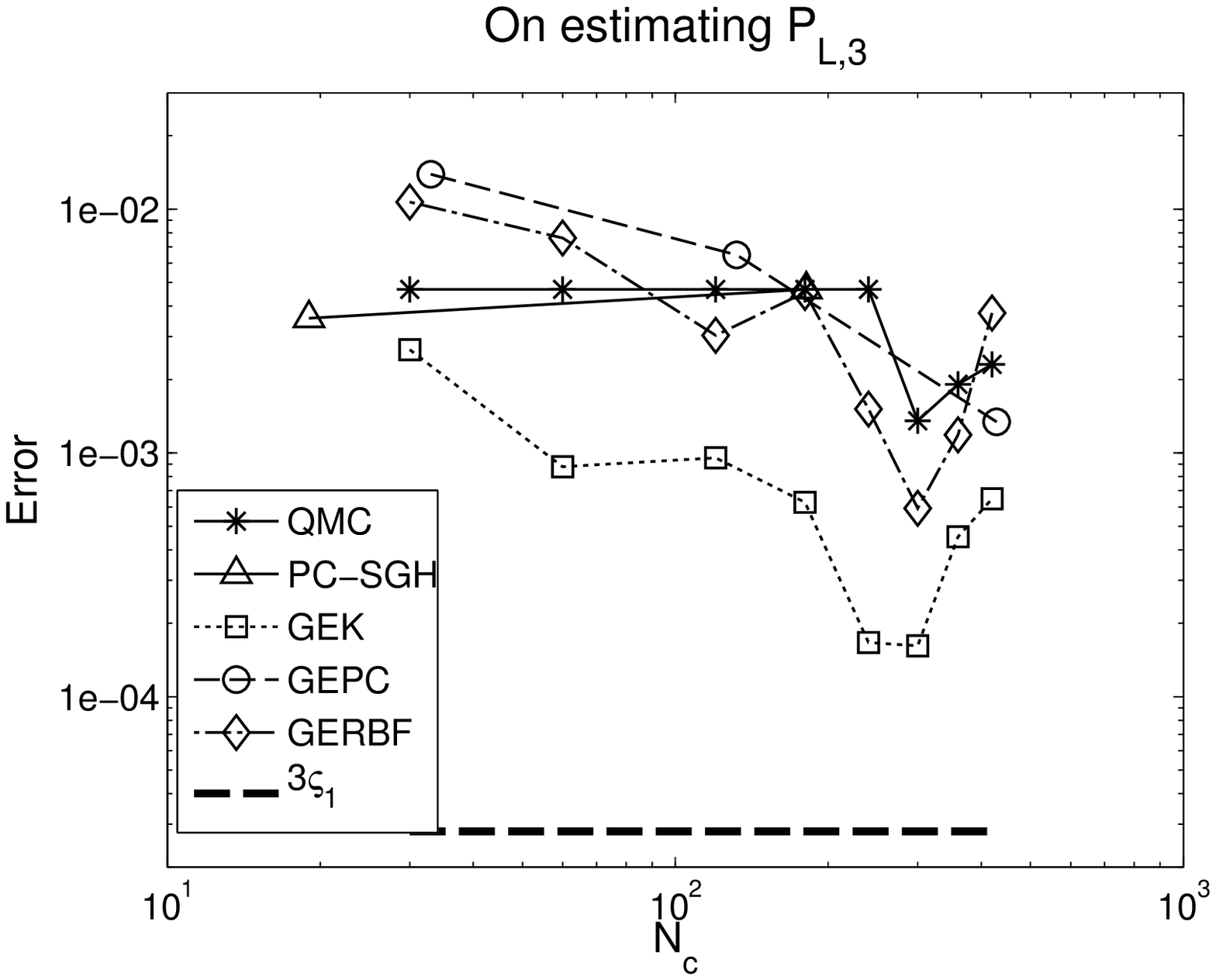}  
\caption{Absolute error in estimating mean, standard deviation (upper row) and exceedance probabilities (lower row) of $C_L$ }
\label{fig:error_CL}
\end{figure}

\begin{figure}[htbp!]
 \centering
 \includegraphics[width=0.43\textwidth]{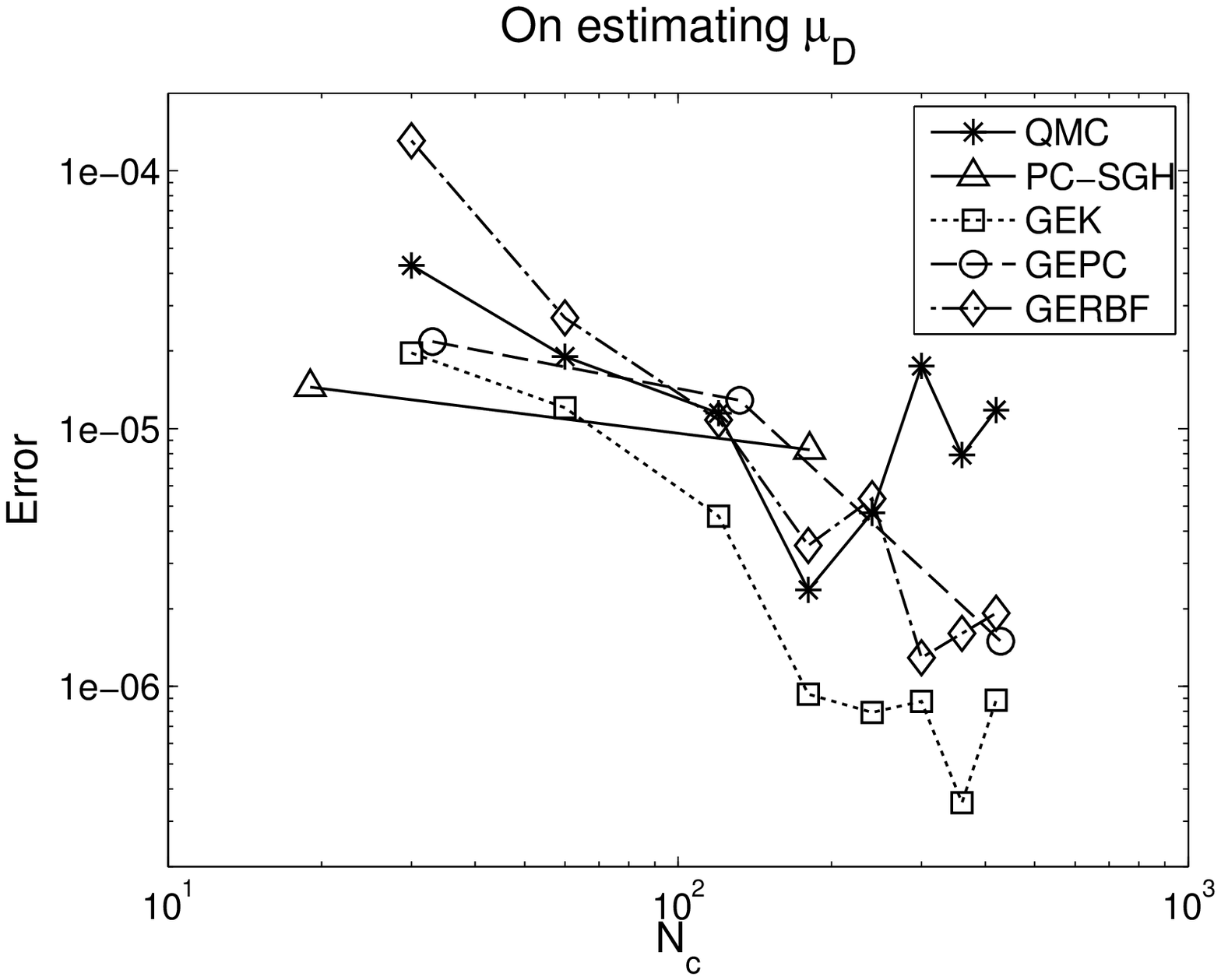}
 \includegraphics[width=0.43\textwidth]{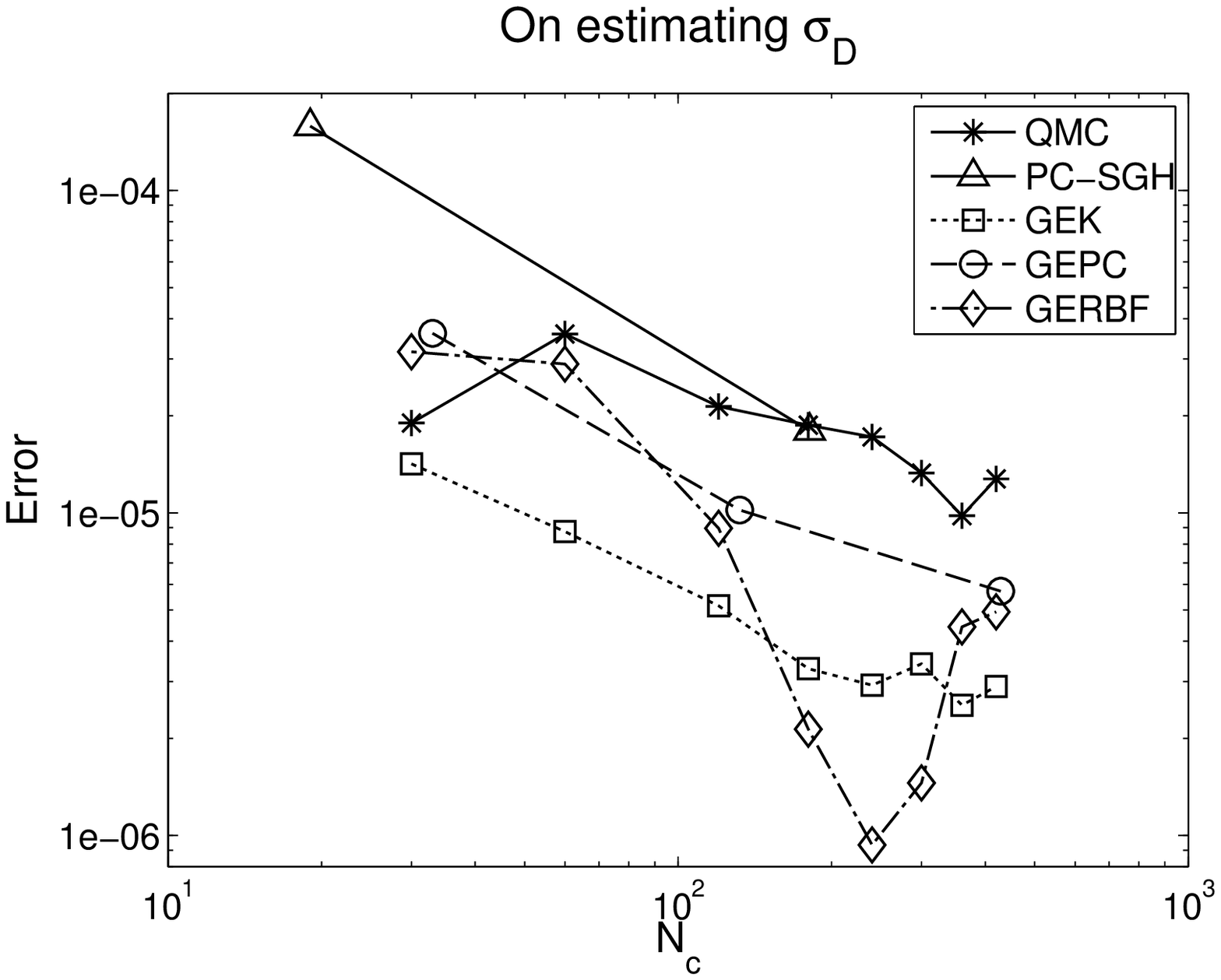} 
 \includegraphics[width=0.43\textwidth]{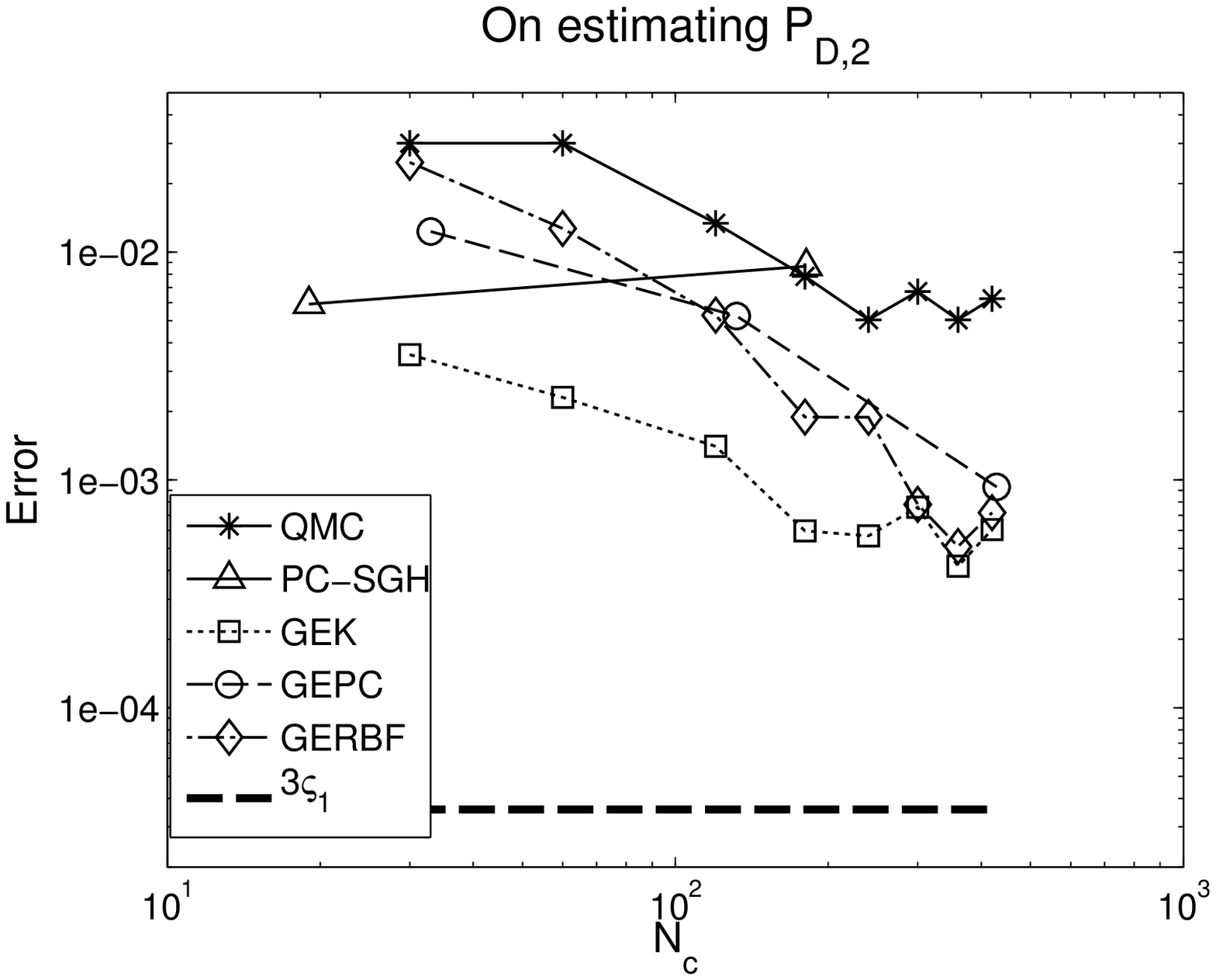}
 \includegraphics[width=0.43\textwidth]{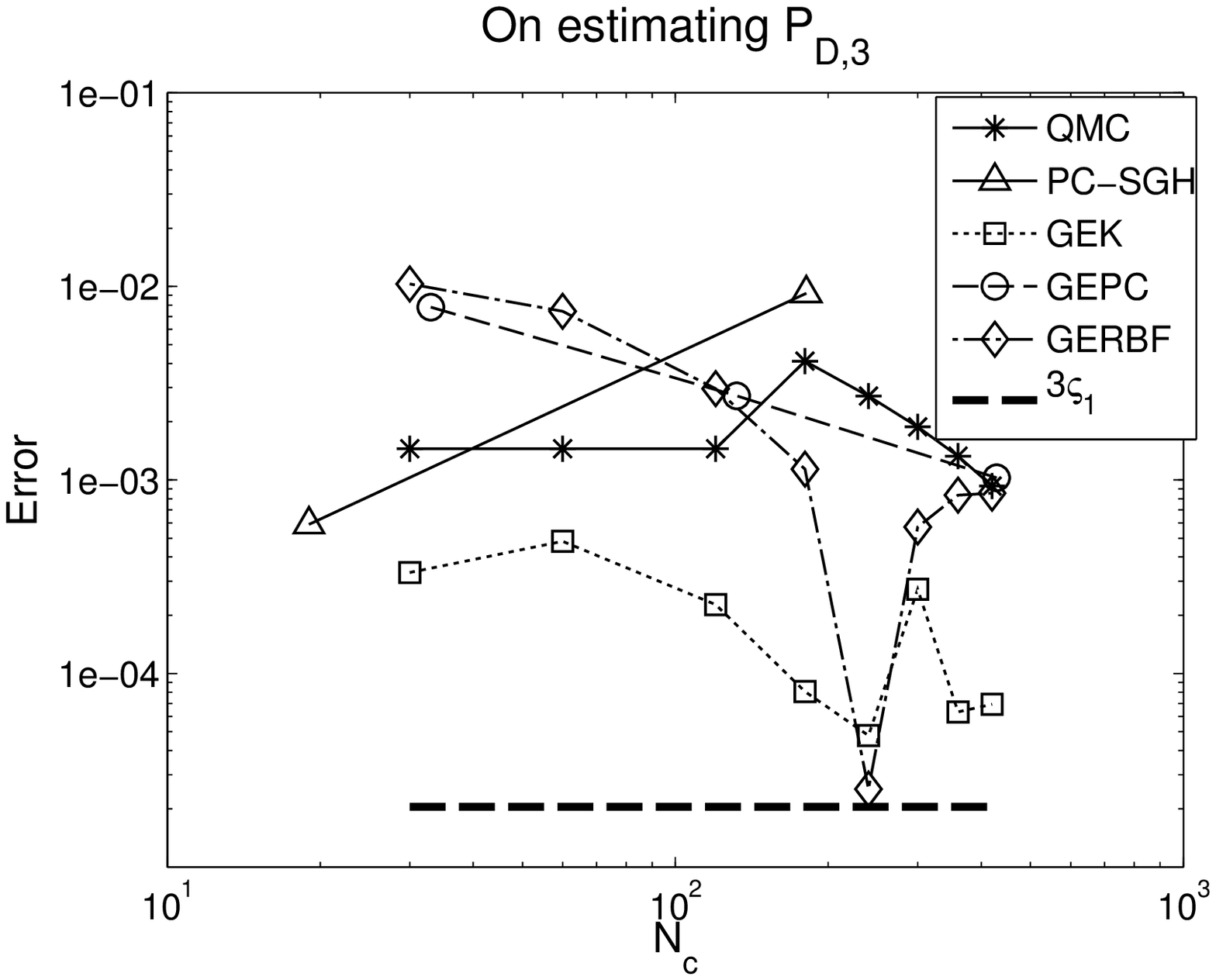}
\caption{Absolute errors in estimating mean, standard deviation (upper row) and exceedance probabilities (lower row) of $C_D$ }
\label{fig:error_CD}
\end{figure}

\begin{figure}[htbp!]
\centering
 \includegraphics[width=0.43\textwidth]{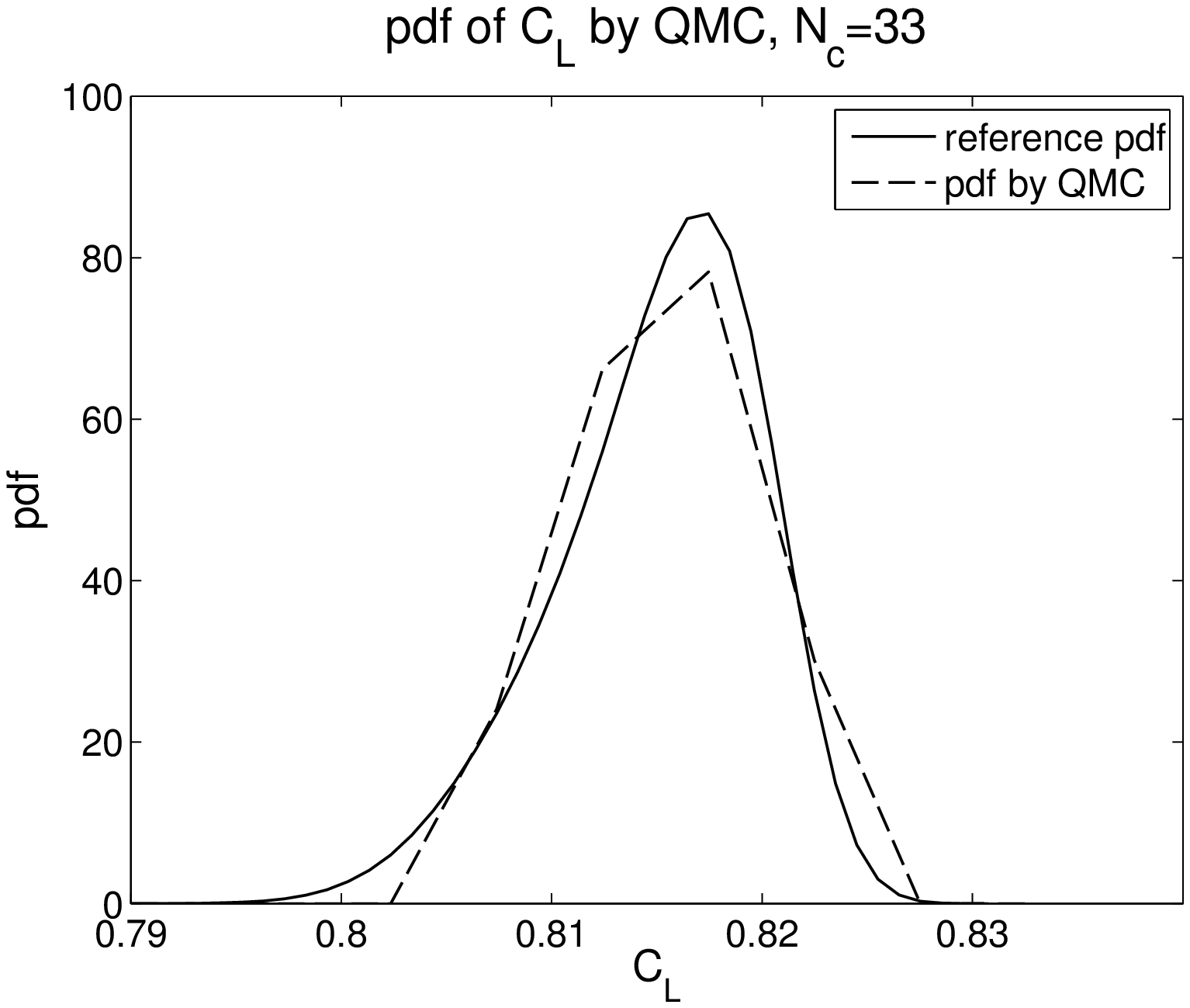} 
 \includegraphics[width=0.43\textwidth]{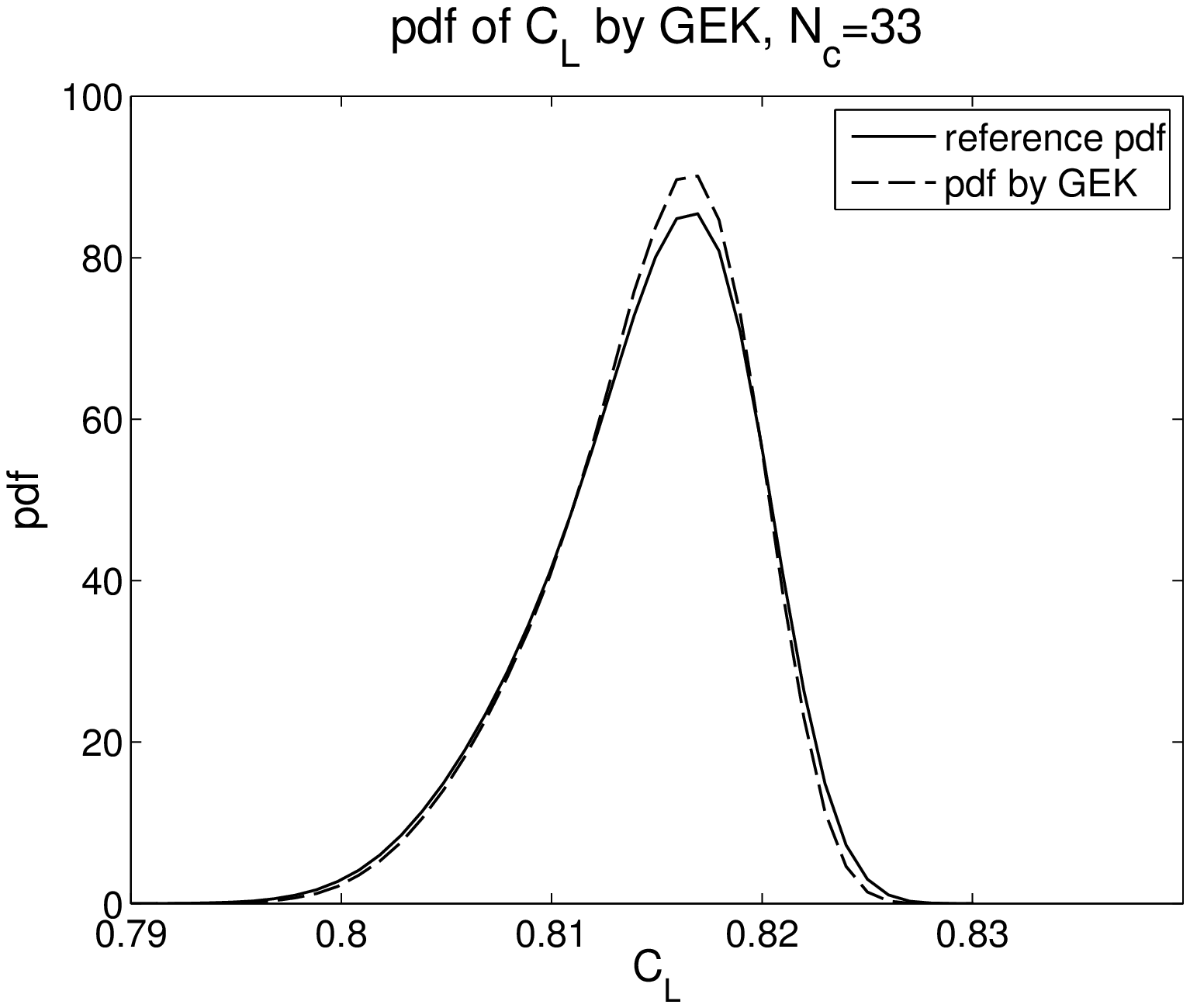} 
 \includegraphics[width=0.43\textwidth]{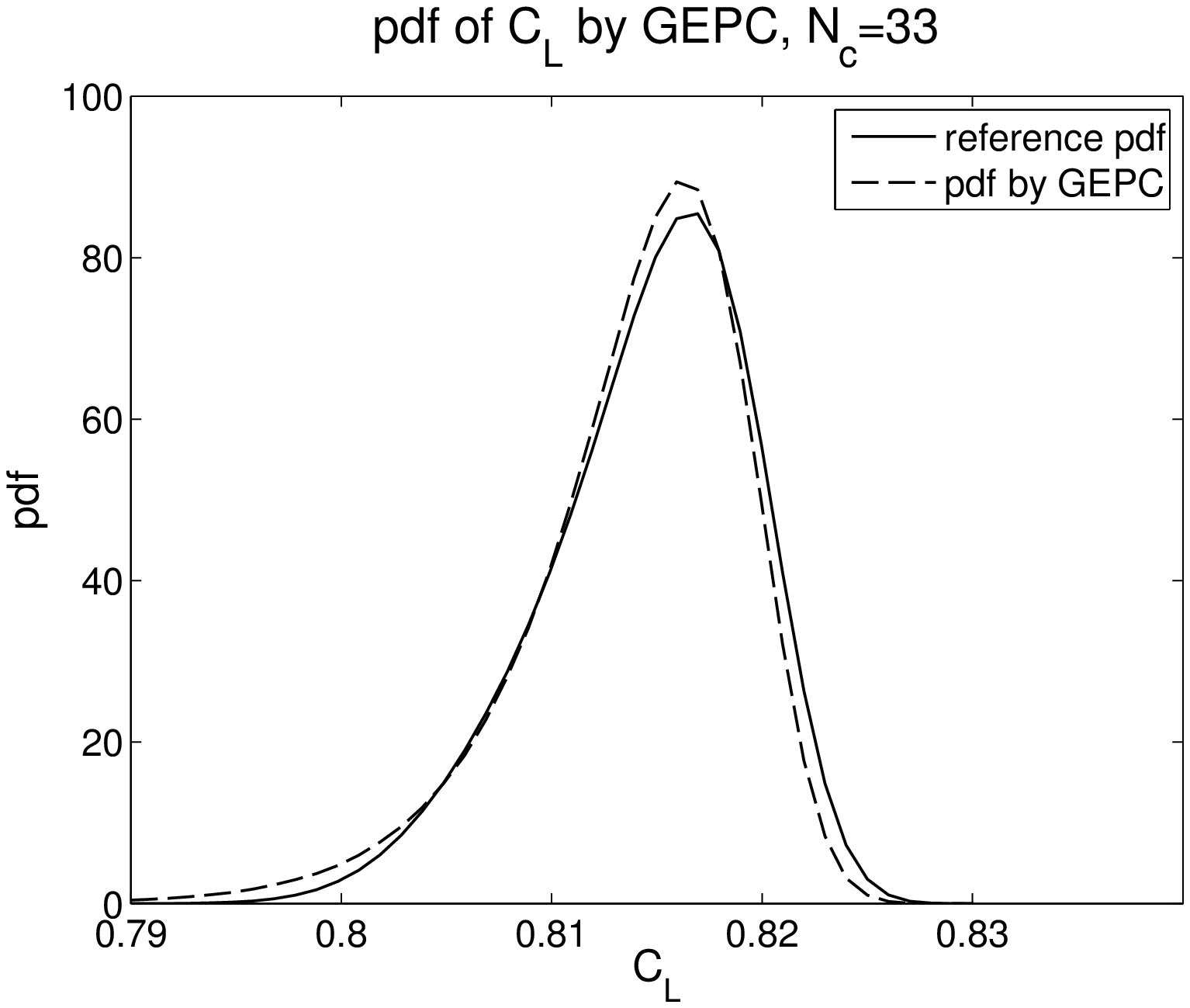} 
 \includegraphics[width=0.43\textwidth]{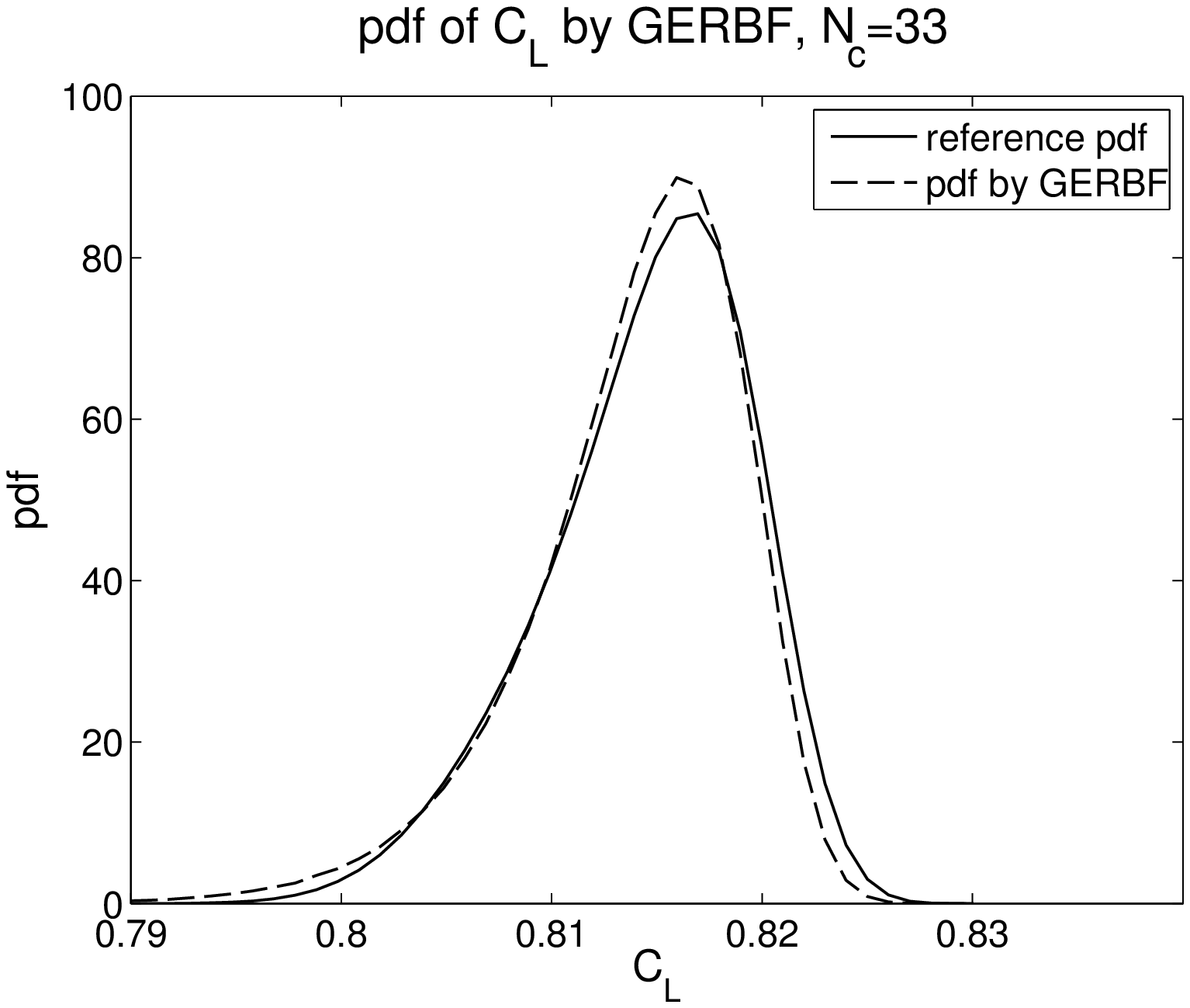} 
\caption{Estimated  pdf (in dash line) of $C_L$ by QMC, GEK, GEPC and GERBF at $N_c=33$.}
\label{fig:pdf_CL}
\end{figure}

\begin{figure}[htbp!]
\centering
  \includegraphics[width=0.43\textwidth]{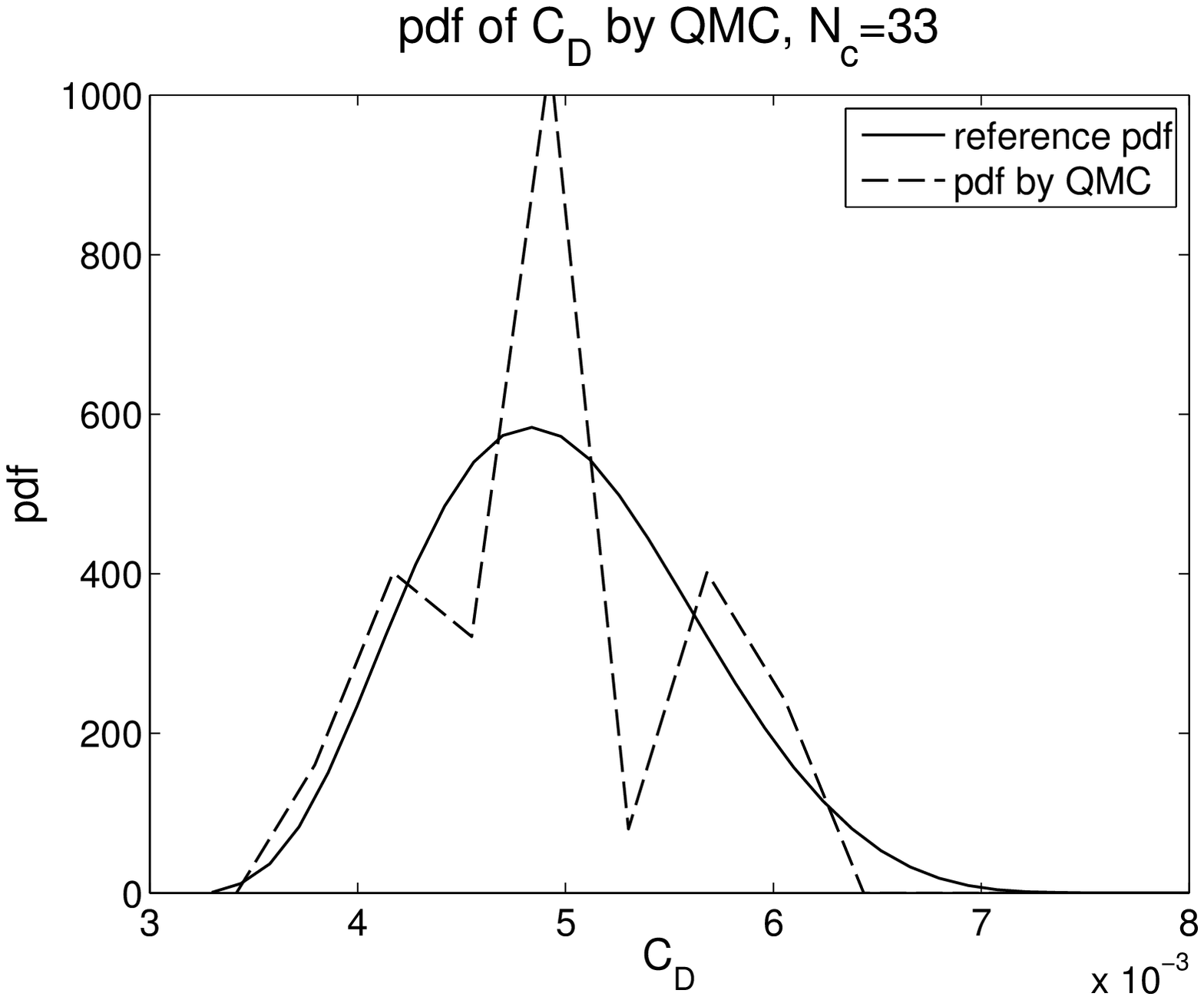} 
 \includegraphics[width=0.43\textwidth]{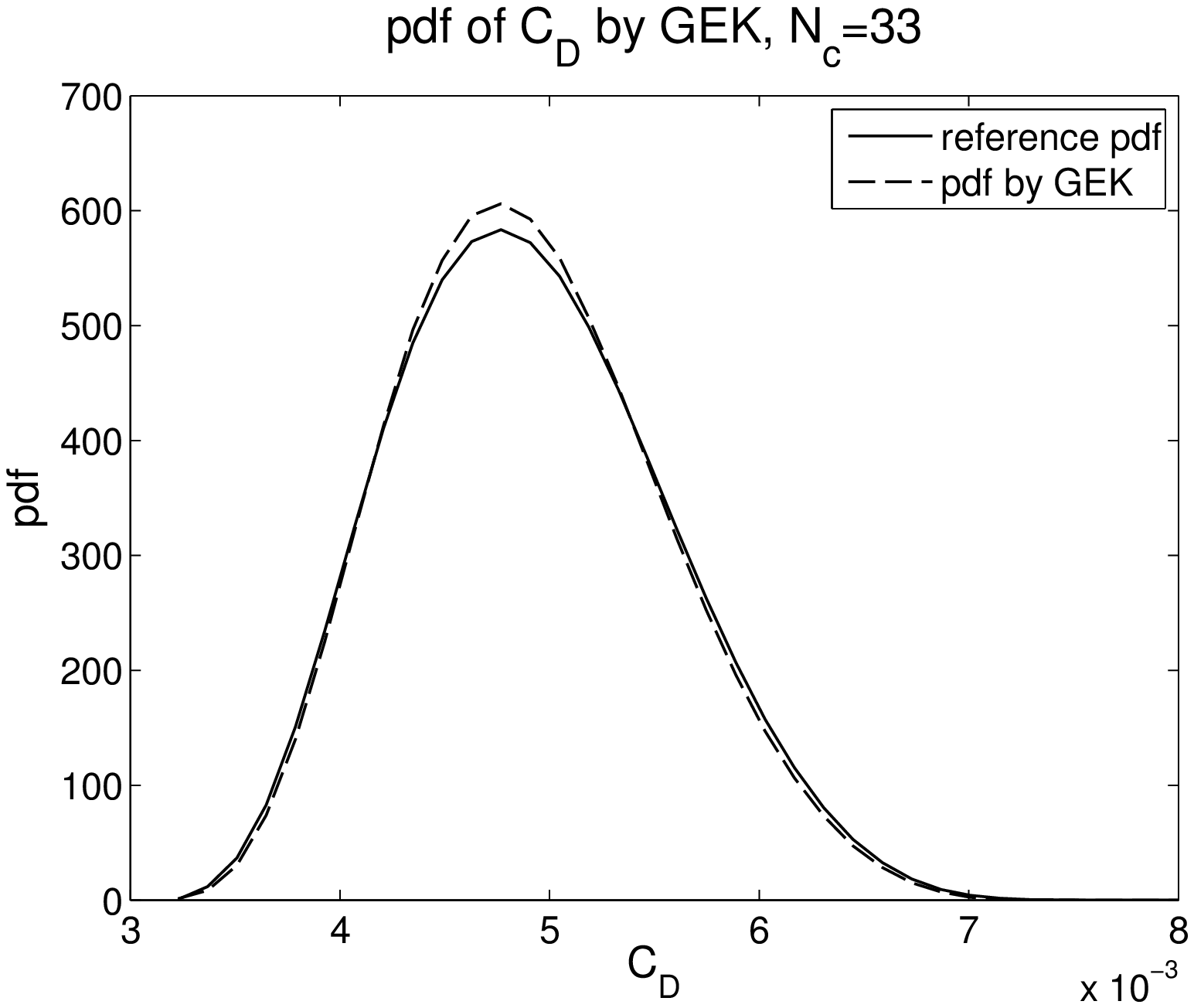} 
 \includegraphics[width=0.43\textwidth]{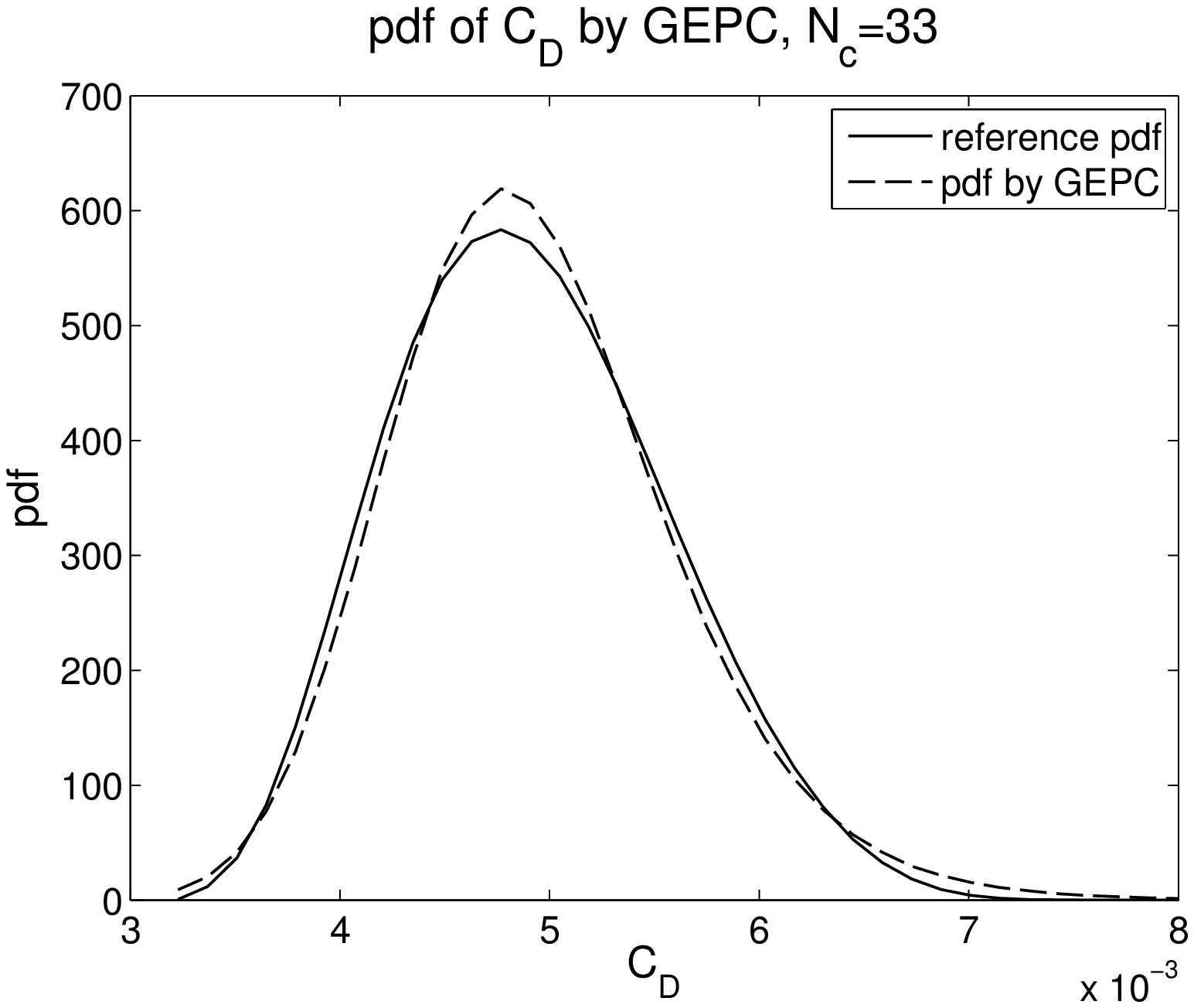} 
  \includegraphics[width=0.43\textwidth]{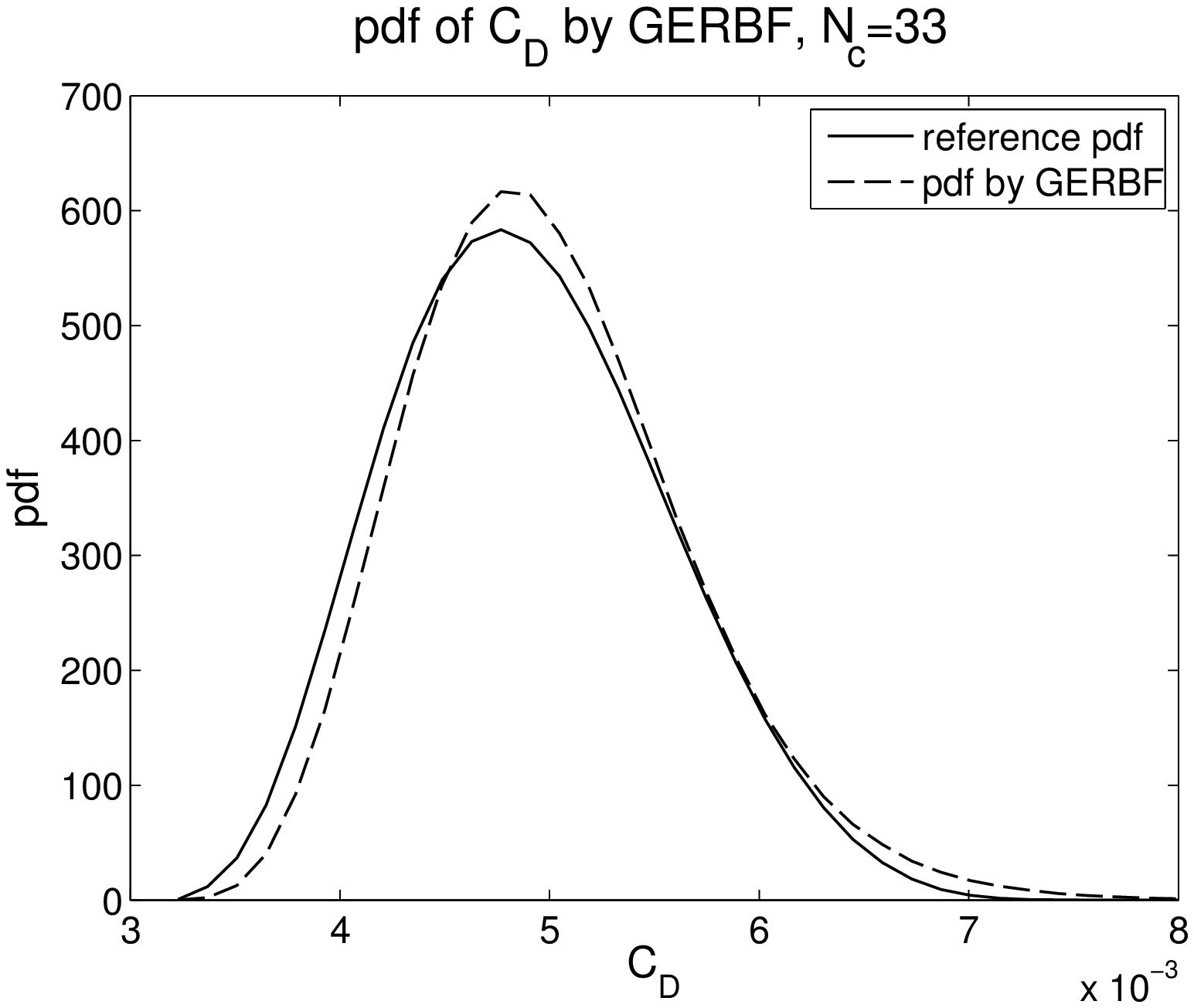} 
\caption{Estimated  pdf (in dash line) of $C_D$ by QMC, GEK, GEPC and GERBF at $N_c=33$.}
\label{fig:pdf_CD}
\end{figure}

The results of the efficiency comparison  are shown in Figures~\ref{fig:error_CL}-\ref{fig:pdf_CD}. Figures \ref{fig:error_CL} and  \ref{fig:error_CD} show  the errors of the five methods in estimating the target statistics of $C_L$ and  $C_D$.   It is observed there that generally the gradient-employing surrogate methods perform better than direct integration methods.   This can be partially ascribed to that the former utilize  more information  with the same computational cost  $N_c$, i.e. they use  $\frac{1+d}{1+s} N_c$  conditions with $s$ the number of SRQs, while a direct integration method uses   $N_c$  conditions. This advantage comes from the cheaper cost of the gradients computed by an adjoint solver when $s$ is smaller than $d$ (in this test case, 2 versus 9),  and the advantage would increase for a larger  $d$. 

The  PC-SGH method  has only two data points  due to the very limited choices of sample number. It is hard to evaluate its error convergence property on only two data points. This   shows a shortcoming of sparse grid quadratures   with a relative high number of variables.  The GEPC method shows constant   convergence in error, but not as efficient as the other two gradient-employing methods. The reason could be that the polynomial surrogate tend to ``overshoot'' in the outskirt of the domain.   But the  PC methods have a merit that they do not need a parameter optimization procedure. A future combination of gradient-employing and the more sophisticated PC methods introduced in \cite{Resmini16} has the potential to enhance the efficiency.

GEK and GERBF are the most efficient methods as far as seen from these results. Beside the cheaper gradients, this  could be attributed to properties of the kernel functions they use  and the effort of tuning parameters. Convergence rate of  inverse multiquadric RBF     was estimated $\C O(e^{- \gamma /h})$ with $h$ the fill distance and $\gamma$ a constant, which translated to a rate in $N$ at $\C O(e^{-\gamma N^{1/d}})$ \cite{Wendland2005}.  We assume the error in the statistics is proportional to that in point-wise approximations,  and find in this 9-variate test case the observed convergence rate is much better than  $\C O(e^{-\gamma N^{1/9}})$,  this hints that the ``effective'' fill distance $h$ reduces faster than $\C O(N^{-1/d})$  due to that some variables are less important than the others (typical for a KLE-based parameterization). This is exploited by anisotropic kernel functions out of the parameter tuning.

GEK seems slightly better than GERBF according to the  results, especially at smaller $N_c$ values. But this might not be generalized to other applications. The advantage of GEK  could come from the possible advantage of cubic spline kernel over inverse multiquadric in this particular case. Different ways of utilizing gradients by GEK (involving 1st and 2nd-order derivatives   and generating a symmetric kernel matrix) and GERBF (involving only 1st order derivative and generating non-symmetric matrix) could also contribute to the difference in their performance.  

At larger $N_c$ values a ``rebound'' of error can be observed in GERBF and GEK results which is caused by numerical instablities due to high condition number of their system matrices, in spite of the stabilizing treatments\footnote{For the GERBF, the truncated SVD as introduced in section \ref{sec:methods}. For the GEK implemented by SMART the system matrix is regularized by adding a  $10^{-13}$  ``nugget'' to the diagonal. They could be insufficient when size of the system matrices become large {\em enough}.} set {\em a priori} in the two methods.  Though there is an \textit{uncertain relation}  (or trade-off principle) exists stating that accuracy and stability cannot be good at the same time \cite{Schaback1995},  more advanced stabilizing techniques,  e.g., pivoted Cholesky decomposition \cite{liu2015} could improve the convergence as sample number is large. 
   
For the mean and standard deviation the reliability of the measured errors are justified by that they are at least by 10 times larger than the corresponding $3\times\varsigma_1$ ($\varsigma_1$ defined in \ref{sec:reference_values}). For exceedance probabilities  the $3\times\varsigma_1$ values are not that small (as indicates by the thick dash lines in the pictures) but are still small enough for most of the measured errors.    

Figures \ref{fig:pdf_CL} and  \ref{fig:pdf_CD} show the probability density functions(pdf) of $C_L$ and  $C_D$ estimated by QMC and the three gradient-employing surrogate methods with  $N_c=33$  (the smallest possible  $N_c$  value for GEPC method), comparing with the reference pdf (computed through  $4 \times 10^6$ QMC samples). There one observes  that for the same computational cost,  the surrogate methods yield  much more accurate pdf's.  This is consistent with their relative performance in estimating the statistics. 

However, the advantage of the gradient-employing surrogate methods should not be taken universal. For example, on problems with few variables or that are mostly undifferentiable the advantage might not exist. If the problem has a very oscillatory topography the advantage could only be manifested with relatively more samples.   

\section{A Navier-Stoke test case of uncertain airfoil geometry}
\label{sec:test case 2}

\begin{figure}[htbp]
  \centering
  \begin{minipage}[b]{0.49\textwidth}
\centering
    \includegraphics[width=\textwidth]{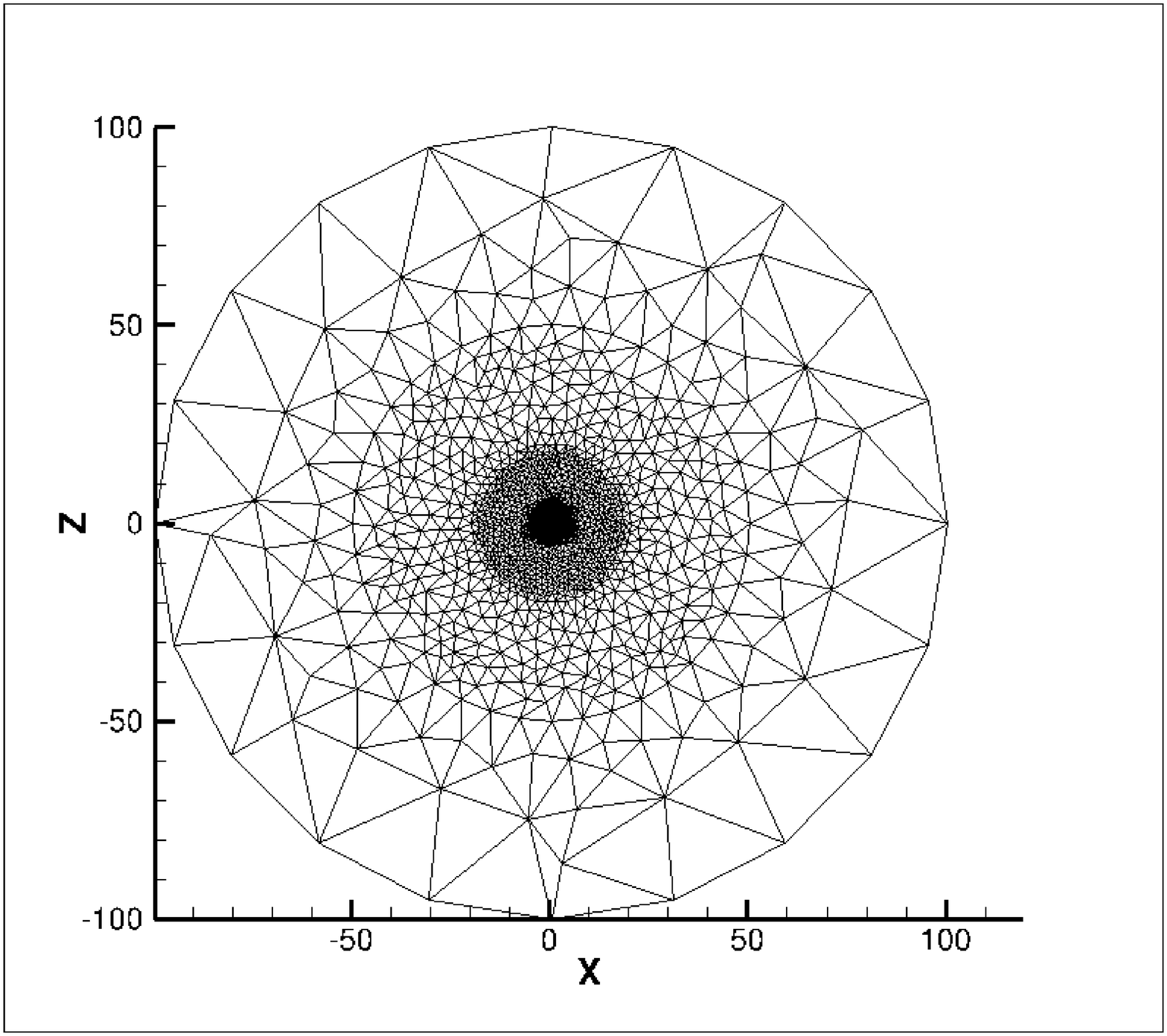}  
  \end{minipage}
  \begin{minipage}[b]{0.49\textwidth}
\centering
    \includegraphics[width=1.0\textwidth]{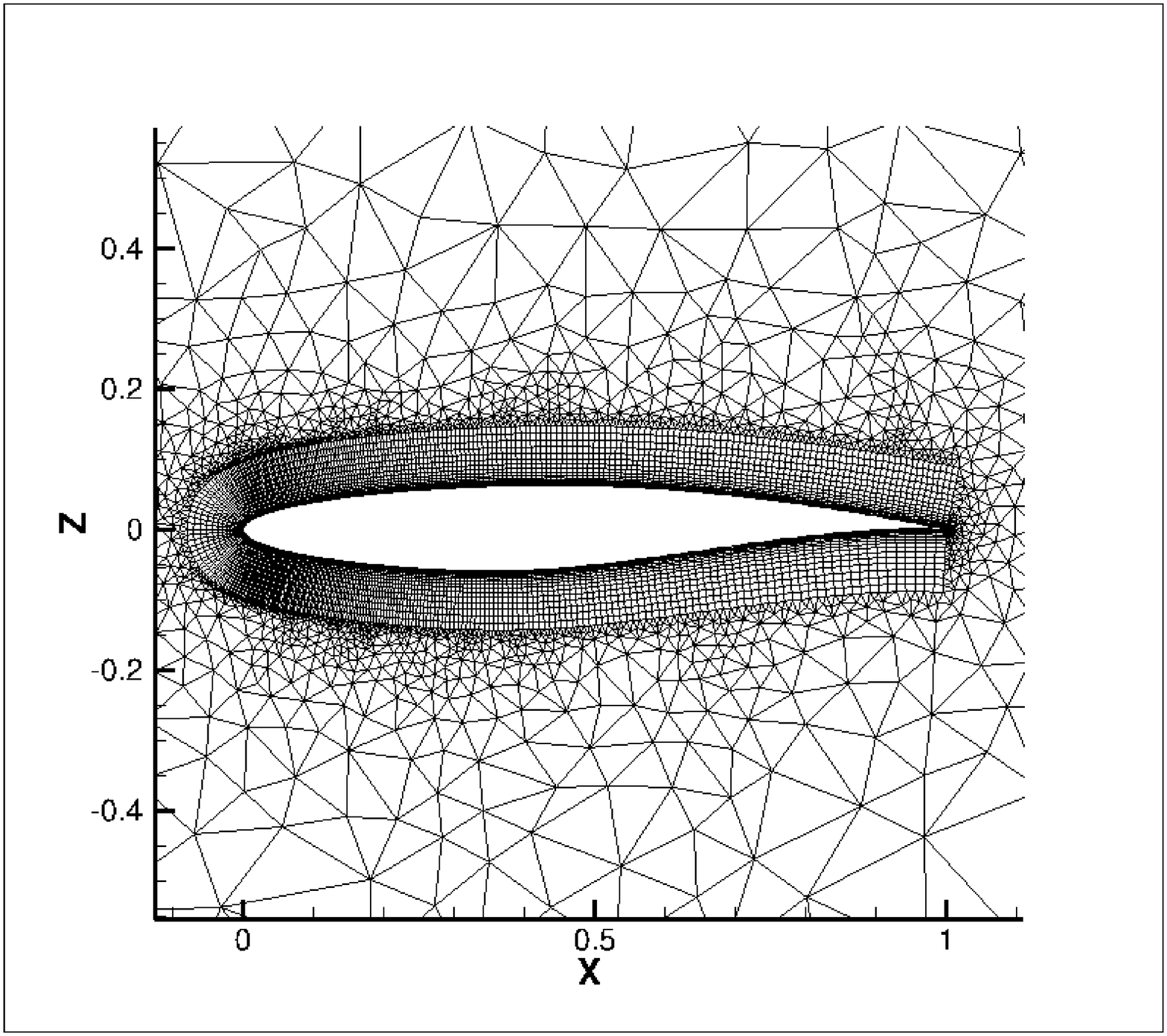}  
  \end{minipage}
\caption{\label{delft_grid} Grid for the RAE2822 airfoil in the RANS test case : the total grid (left) and zoom around the airfoil (right).}
\end{figure}

The second test case is also based on a RAE2822 airfoil, but here we use a Reynolds averaged Navier--Stokes(RANS) solver in the TAU package\cite{TAU2}, opting for a turbulence model SAO, with a nominal Mach number 0.729 and angle of attack(AoA) 2.31 degree.  The domain is discretized by a much finer grid in which the airfoil has 380 surface nodes, as shown in  Figure~\ref{delft_grid}. 

\begin{figure}[htbp!]
\centering
\includegraphics[width=0.48\textwidth,height=9cm]{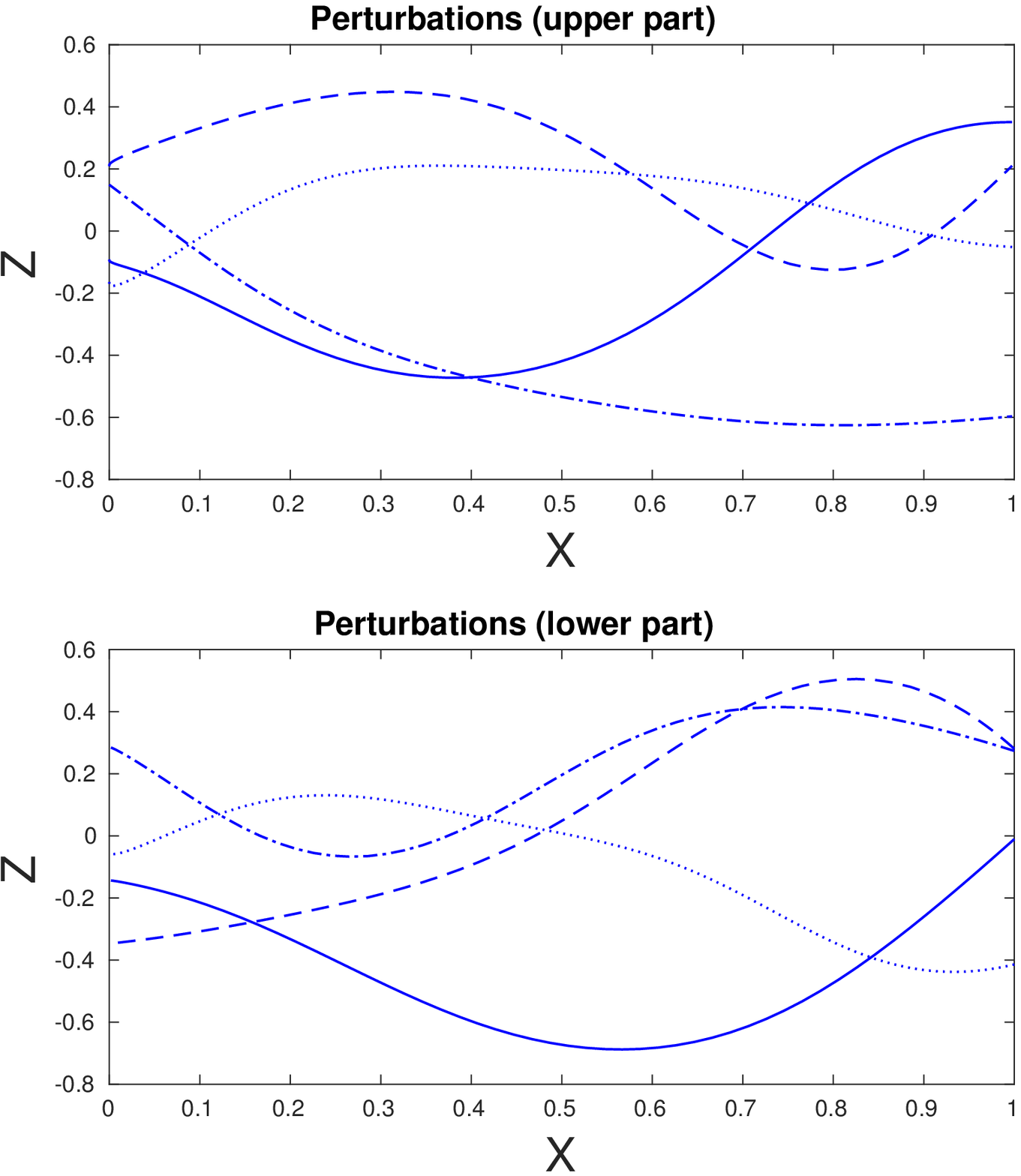}
\includegraphics[width=0.48\textwidth,height=9cm]{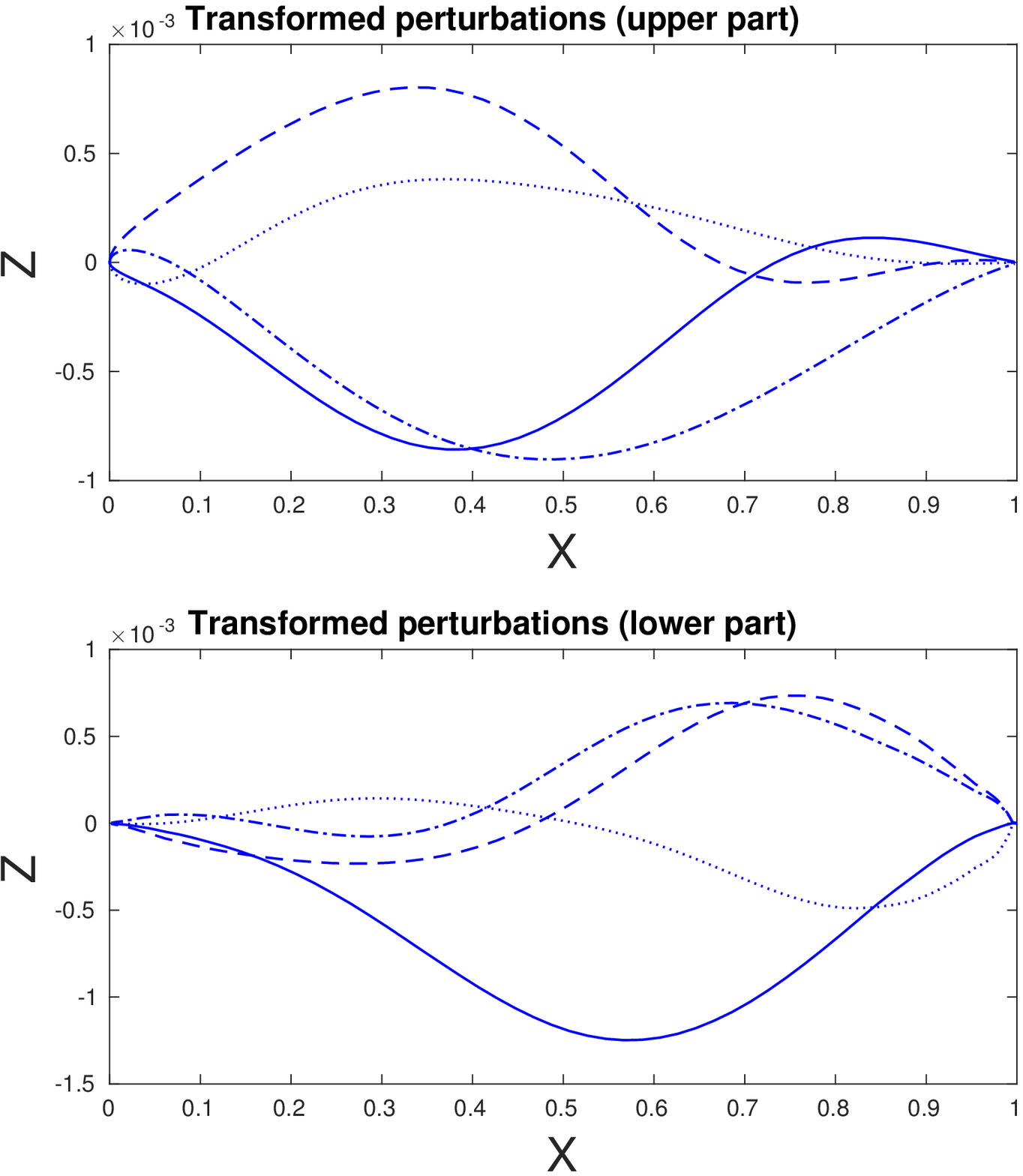} 
\caption{\label{Gaussian_field} Four realizations of the Gaussian random field $\gaussfield$ (left) and transformed field $\vek R$  (right).}
\end{figure}

\begin{figure}[htbp]
    \includegraphics[width=0.7\textwidth]{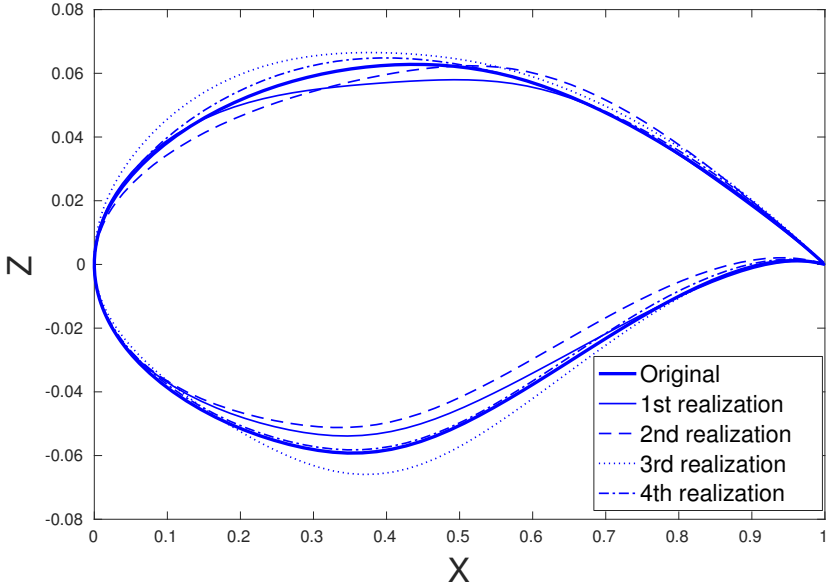}  
\caption{\label{resulpertshapes_2} Four realizations of perturbed RAE2822 airfoil geometry, 10 times exaggerated. }
\end{figure}

The sources of uncertainty are random perturbation in the airfoil geometry and operational parameters Mach number and AoA. The geometry perturbation is modeled by a zero-mean random field through a Karhunen-Lo\`{e}ve expansions (KLE) using the same methodology as introduced in the first test case, but with a different setting. Here the upper and lower surfaces of the airfoil are treated as two separate random fields as their correlation in geometric variations is assumed weak in this test case.  The standard deviation $\sigma $ in (\ref{cov}) is uniformly set to $1$ in the fields. After the standard Gaussian random fields are approximated through the KLE in (\ref{KLE}), they are transformed to perturbation fields $\vek R$ being of Beta distribution\footnote { The  probability density function is $P(x)=\frac{(x_{max}-x)^{(\beta-1)}(x-x_{min})^{(\alpha-1)}}{B(\alpha, \beta)}$ with $[x_{min}, x_{max}]$ support of the distribution, $B$ Beta function and $\alpha = \beta = 4$. }  within support $[-\delta(\vek g), \delta(\vek g) ]$ for the boundedness. $\delta(\vek g)$ is set to $1.5\%$ of the airfoil thickness  at the node $\vek g$. We keep all the eigenvalues larger than $10^{-7}$ in the truncated  KLE approximations, this leads to a parametrisation in 24 variables for the geometric perturbation only.  Figure~\ref{Gaussian_field} display four realizations of the $\gaussfield$ and $\vek R$ respectively.  Examples of perturbed  geometries  $\widetilde {\mat G}$  are shown in Figure~\ref{resulpertshapes_2}. 

The perturbations in Mach number and AoA are of the same Beta distribution, but with a support within $\pm 2\%$ of the nominal values, adding up to a total number of variables as 26 for the test case. 

Three UQ methods,  quasi-Monte Carlo quadrature, gradient-enhanced kriging and plain kriging, are applied to the test case  and their efficiency compared in estimating four statistics (mean, standard deviation, skewness and kurtosis) of the coefficient of lift ($C_L$).

The accuracy of the estimates is judged by comparing to reference statistics which are based on 10000 QMC samples.  

The computational cost is again measured in terms of ``compensated evaluation number'' $N_c$ as introduced in the first test case. The only difference is that for the gradient-employing method (GEK) $N_c$ is set to $2N$ since here we only handle one system response quantity (SRQ) $C_L$.  

\subsection{Numerical results and discussion }\label{sec:results_2}
\begin{figure}[htbp!]
  \centering
 \includegraphics[width=0.43\textwidth]{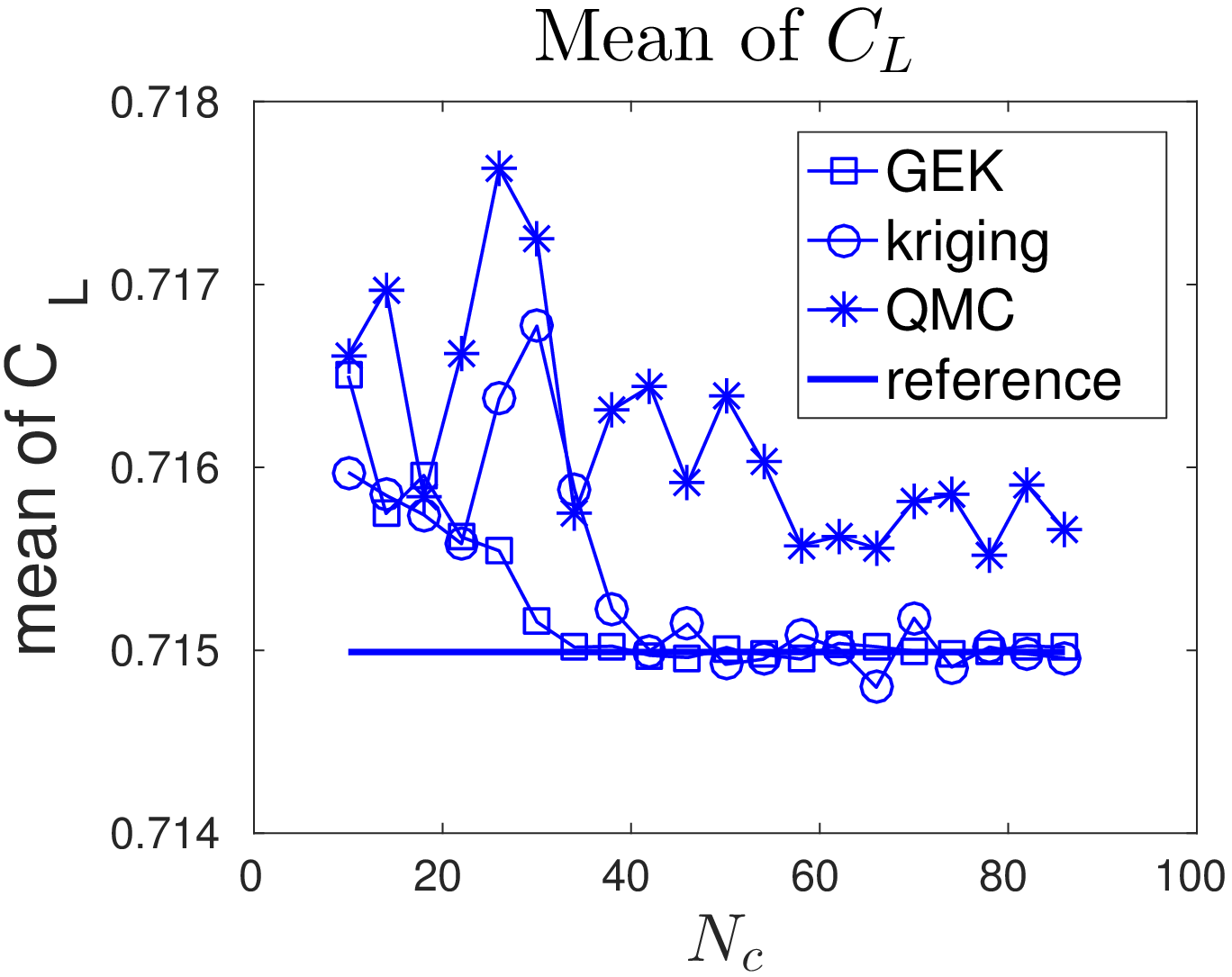}
 \includegraphics[width=0.43\textwidth]{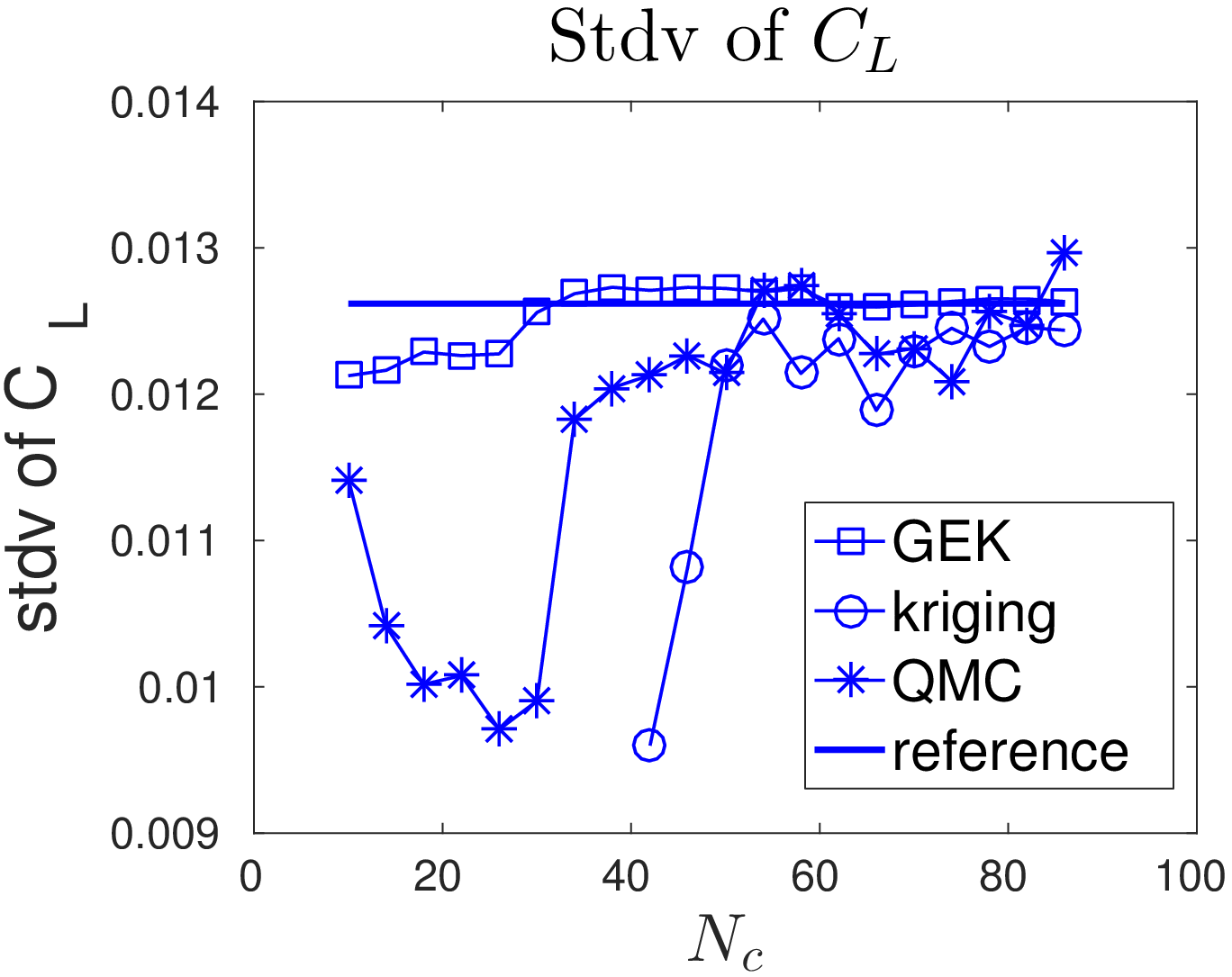}  
 \includegraphics[width=0.43\textwidth]{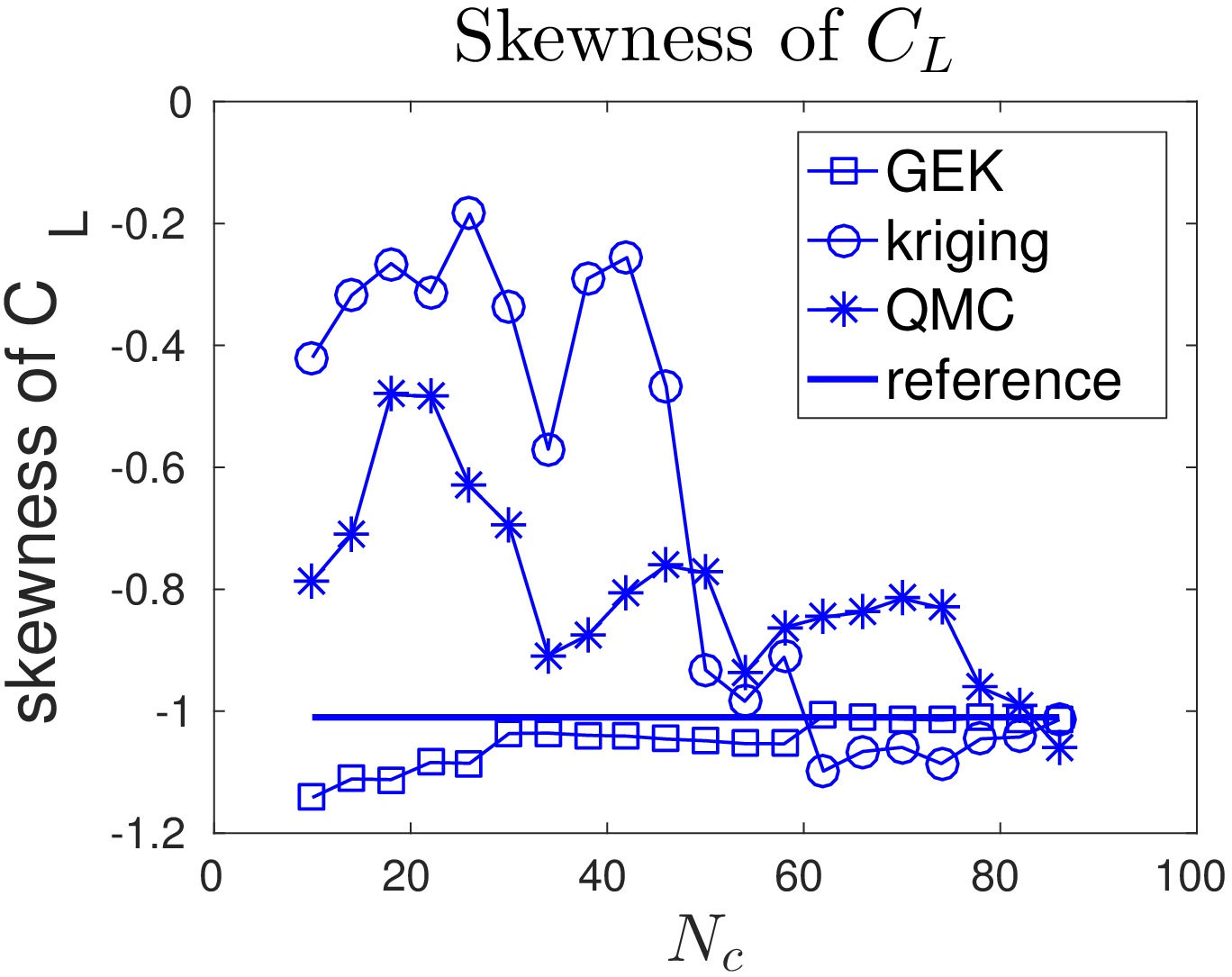}
 \includegraphics[width=0.43\textwidth]{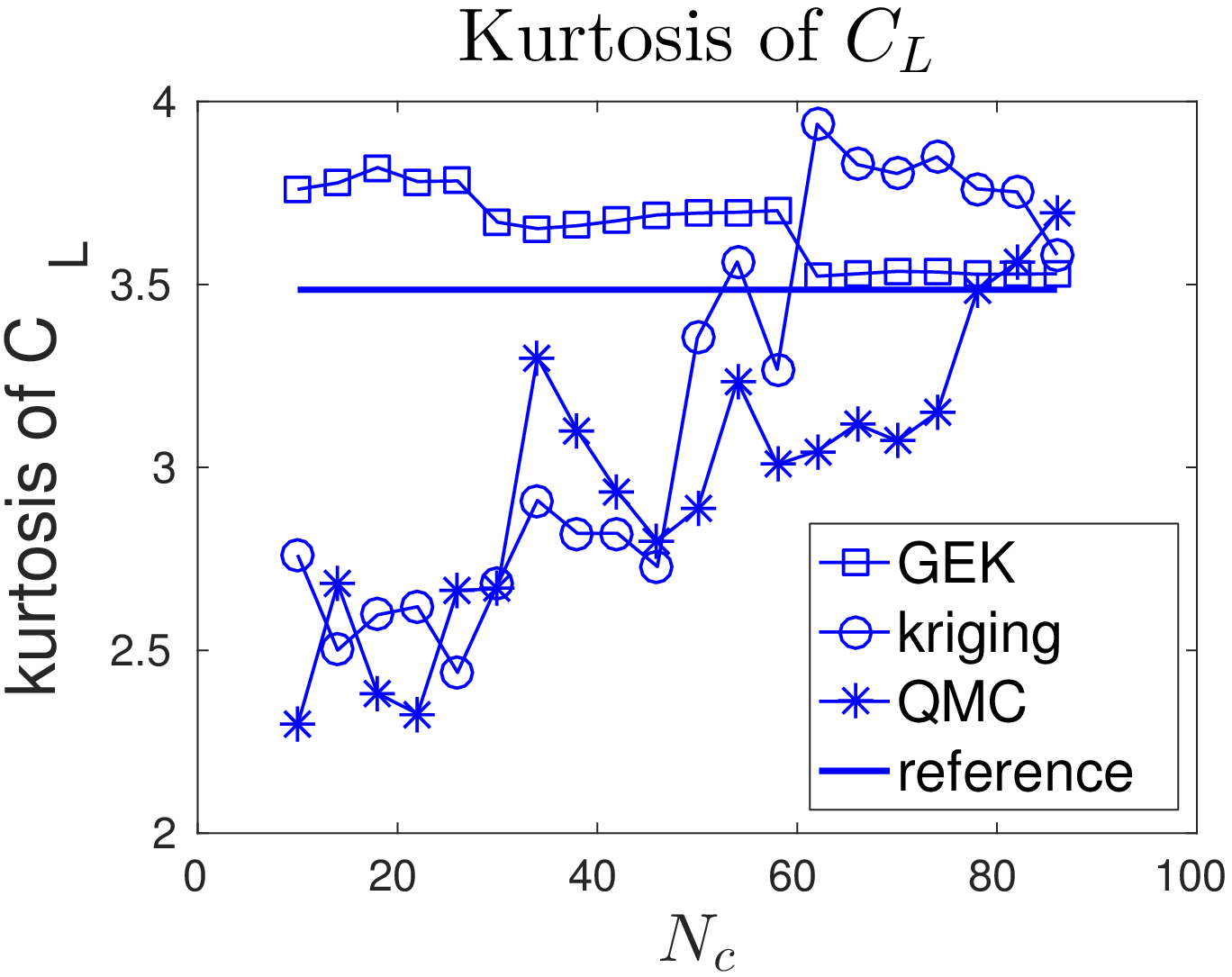}  
\caption{ Estimates of mean, standard deviation, skewness and kurtosis of $C_L$, contrasted to the reference.}
\label{fig:error_CL_2}
\end{figure}

For this RANS test case our reference statistics at hands are not as accurate as those in the Euler one due to the much fewer CFD evaluations affordable. Therefore displaying error convergence in logarithmic scale is not fully supported by the precision. So we are content with observing the approaching of the estimated statistics to the references in a real scale.  In Figures~\ref{fig:error_CL_2} we contrast the estimates of the statistics to the reference statistics. As expected, generally the estimates approach to the references along the increase of the cost measure $N_c$. It is observed that the gradient-enhanced kriging (GEK) outperforms plain kriging and quasi-Monte Carlo quadrature(QMC) in the estimation of all statistics. As in the first test case this can be explained by the more information utilized by the former method at the same computational cost. Especially the contrast of GEK and plain kriging highlights the efficiency gain brought by the cheaper gradient information. 

 
\section{Conclusion} \label{sec:conclusion}
 
This paper compares efficiency of various methods in quantifying geometry-induced aerodynamic uncertainties. An Euler and a RANS test cases  based on RAE2822 airfoil are set up by perturbing the   geometry by  random fields which are approximated and parametrized through Karhunen-Lo\`{e}ve expansions. UQ methods including Quasi-Monte Carlo  quadrature, polynomial chaos with coefficients determined by sparse grid quadrature, gradient-enhanced radial basis functions, gradient-enhanced polynomial chaos, gradient-enhanced kriging and plain kriging are applied to the test cases and are compared in their efficiency in estimating some statistics and probability distribution of the uncertain lift and drag coefficients. The results show that gradient-employing surrogate methods are more efficient. The advantage is due to the cheaper gradients obtained by using adjoint solver, and is expected to increase with an increasing $d$.   
 
\section*{Acknowledgements}
This research has been conducted within the project MUNA under the framework of the German Luftfahrtforschungsprogramm funded by the Ministry of Economics (BMWi). 

Additionally, a part of this work was done by A. Litvinenko during his stay at King Abdullah University of Science and Technology.

We are grateful to the anonymous reviewers for their diligence and insights which help greatly to improve this paper. A special gratitude is extended to Dr. Stefan G\"{o}rtz in DLR for his invaluable advices during the modification.

\FloatBarrier
 
\bibliographystyle{siam}      
\bibliography{referenceMUNA}{}   

\begin{thebibliography}{10}

\bibitem{adler2007}
{\sc R.~J. Adler and J.~E. Taylor}, {\em Random Fields and Geometry},
  Springer-Verlag, Berlin, 2007.

\bibitem{Barthelmann2000}
{\sc V.~Barthelmann, E.~Novak, and K.~Ritter}, {\em High dimensional polynomial
  interpolation on sparse grids}, Advances in Computational Mathematics, 12
  (2000), pp.~273--288.

\bibitem{Bijl1}
{\sc H.~Bijl, D.~Lucor, S.~Mishra, and C.~Schwab}, {\em Uncertainty
  quantification in computational fluid dynamics}, Lecture Notes in
  Computational Science and Engineering, Springer, 2013.

\bibitem{bompard2010}
{\sc M.~Bompard, J.~Peter, and J.-A. D\'{e}sid\'{e}ri}, {\em Surrogate models
  based on function and derivative values for aerodynamic global optimization},
  in Fifth European Conference on Computational Fluid Dynamics, ECCOMAS CFD
  2010, Lisbon, Portugal, 2010.

\bibitem{brezillon2005}
{\sc J.~Brezillon and R.~Dwight}, {\em Discrete adjoint of the navier-stokes
  equations for aerodynamic shape optimization}, in EUROGEN 2005 - Sixth
  Conference on Evolutionary and Deterministic Methods for Design, Optimization
  and Control with Applications to Industrial and Societal Problems, Munich,
  Germany, 2005.

\bibitem{buhmann2000}
{\sc M.~D. Buhmann}, {\em Radial basis functions}, Acta Numerica, 9 (2000),
  pp.~1--38.

\bibitem{caf1998}
{\sc R.~E. Caflisch}, {\em Monte {C}arlo and quasi-{M}onte {C}arlo methods},
  Acta Numerica, 7 (1998), pp.~1--49.

\bibitem{Chen2011}
{\sc S.~Chen and W.~Chen}, {\em A new level-set based approach to shape and
  topology optimization under geometric uncertainty}, Structural and
  Multidisciplinary Optimization, 44 (2011), pp.~1--18.

\bibitem{chung2002}
{\sc H.-S. Chung and J.~J. Alonso}, {\em Using gradients to construct cokriging
  approximation models for high-dimensional design optimization problems}, AIAA
  paper, 317 (2002), pp.~14--17.

\bibitem{dolgov2014computation}
{\sc S.~Dolgov, B.~N. Khoromskij, A.~Litvinenko, and H.~G. Matthies}, {\em
  Computation of the response surface in the tensor train data format}, arXiv
  preprint arXiv:1406.2816,  (2014).

\bibitem{Evans07}
{\sc T.~Evans, P.~Tattersall, and J.~Doherty}, {\em Identification and
  quantification of uncertainty sources in aircraft-related cfd computations ?
  an industrial perspective.}, in Proceedings of RTO-AVT-147 Symposium.
  Athens., Greece, 2007. Paper 6., 2007.

\bibitem{gian2006}
{\sc K.~C. Giannakoglou, D.~I. Papadimitriou, and I.~C. Kampolis}, {\em
  Aerodynamic shape design using evolutionary algorithms and new
  gradient-assisted metamodels}, Computer Methods in Applied Mechanics and
  Engineering, 195 (2006), pp.~6312--6329.

\bibitem{Giraldi2014}
{\sc L.~Giraldi, A.~Litvinenko, D.~Liu, H.~G. Matthies, and A.~Nouy}, {\em To
  be or not to be intrusive? the solution of parametric and stochastic
  equations---the "plain vanilla" galerkin case}, SIAM Journal on Scientific
  Computing, 36 (2014), pp.~A2720--A2744.

\bibitem{giunta2004}
{\sc A.~A. Giunta, M.~S. Eldred, and J.~P. Castro}, {\em Uncertainty
  quantification using response surface approximation}, in 9th ASCE Specialty
  Conference on Probabolistic Mechanics and Structural Reliability,
  Albuquerque, New Mexico, USA, 2004.

\bibitem{han2012}
{\sc Z.~H. Han, S.~G\"{o}rtz, and R.~Zimmermann}, {\em Improving
  variable-fidelity surrogate modeling via gradient-enhanced kriging and a
  generalized hybrid bridge function}, Journal of Aerospace Science and
  Technology,  (2012).

\bibitem{TAU2}
{\sc R.~Heinrich, R.~Dwight, M.~Widhalm, and A.~Raichle}, {\em Algorithmic
  developments in {TAU}}, in {MEGAFLOW} - {N}umerical Flow Simulation for
  Aircraft Design, N.~Kroll and J.~Fassbender, eds., vol.~89, Springer, 2005,
  pp.~93--108.

\bibitem{hosder2007}
{\sc S.~Hosder, R.~W. Walters, and M.~Balch}, {\em Efficient sampling for
  non-intrusive polynomial chaos applications with multiple uncertain input
  variables}, in Proceedings of the 48th AIAA/ASME/ASCE/AHS/ASC Structures,
  Structural Dynamics, and Materials Conference, number AIAA-2007-1939,
  Honolulu, HI, vol.~125, 2007.

\bibitem{hosder2006}
{\sc S.~Hosder, R.~W. Walters, and R.~Perez}, {\em A non-intrusive polynomial
  chaos method for uncertainty propagation in cfd simulations}, in Proceedings
  of the 44th AIAA Aerospace Sciences Meeting, vol.~14, 2006, pp.~10649--10667.

\bibitem{joe2008}
{\sc S.~Joe and F.~Kuo}, {\em Constructing sobol sequences with better
  two-dimensional projections}, SIAM Journal on Scientific Computing, 30
  (2008), pp.~2635--2654.

\bibitem{khoromskij2009application}
{\sc B.~N. Khoromskij, A.~Litvinenko, and H.~G. Matthies}, {\em Application of
  hierarchical matrices for computing the karhunen--lo{\`e}ve expansion},
  Computing, 84 (2009), pp.~49--67.

\bibitem{Kim2006}
{\sc N.~H. Kim, H.~Wang, and N.~V. Queipo}, {\em Efficient shape optimization
  under uncertainty using polynomial chaos expansions and local sensitivities},
  AIAA journal, 44 (2006), pp.~1112--1116.

\bibitem{Lucor08}
{\sc J.~Ko, D.~Lucor, and P.~Sagaut}, {\em Sensitivity of two-dimensional
  spatially developing mixing layers with respect to uncertain inflow
  conditions}, Physics of Fluids, 20 (2008).

\bibitem{le2010spectral}
{\sc O.~Le~Ma{\^\i}tre and O.~M. Knio}, {\em Spectral methods for uncertainty
  quantification: with applications to computational fluid dynamics}, 2010.

\bibitem{LitvSampling13}
{\sc A.~Litvinenko, H.~Matthies, and T.~A. El-Moselhy}, {\em Sampling and
  low-rank tensor approximation of the response surface}, in Monte Carlo and
  Quasi-Monte Carlo Methods 2012, J.~Dick, F.~Y. Kuo, G.~W. Peters, and I.~H.
  Sloan, eds., vol.~65 of Springer Proceedings in Mathematics $\&$ Statistics,
  Springer Berlin Heidelberg, 2013, pp.~535--551.

\bibitem{Litvin_Creta10}
{\sc A.~Litvinenko and H.~G. Matthies}, {\em Low-rank data format for
  uncertainty quantification}, in International Conference on Stochastic
  Modeling Techniques and Data Analysis Proceedings. Editor: Chr. H. Skiadas,
  Chania, Greece, 2010, pp.~477--484.

\bibitem{Litvin11PAMM}
{\sc A.~Litvinenko and H.~G. Matthies}, {\em Uncertainties quantification and
  data compression in numerical aerodynamics}, PAMM, 11 (2011), pp.~877--878.

\bibitem{liu2012b}
{\sc D.~Liu}, {\em {A Best Practice Guide: Efficient Quantification of
  Aerodynamic Uncertainties}}, tech. report, Center for Computer Applications
  in Aerospace Science and Engineering, German Aerospace Center (DLR), 2012.
\newblock IB 124-2012/2.

\bibitem{liu2014}
{\sc D.~Liu and S.~G\"{o}rtz}, {\em Efficient quantification of aerodynamic
  uncertainty due to random geometry perturbations}, in New Results in
  Numerical and Experimental Fluid Mechanics IX, A.~Dillmann et~al., eds.,
  Springer International Publishing, 2014, pp.~65--73.

\bibitem{liu2015}
{\sc D.~Liu and H.~G. Matthies}, {\em Pivoted cholesky decomposition by cross
  approximation for efficient solution of kernel systems}.

\bibitem{Bijl}
{\sc G.~J.~A. Loeven and H.~Bijl}, {\em Airfoil analysis with uncertain
  geometry using the probabilistic collocation method}, in 48th
  AIAA/ASME/ASCE/AHS/ASC Structures, Structural Dynamics, and Materials
  Conference, AIAA-2008-2070, 2008.

\bibitem{Loeven2007}
{\sc G.~J.~A. Loeven, J.~A.~S. Witteveen, and H.~Bijl}, {\em A probabilistic
  radial basis function approach for uncertainty quantification}, in
  Proceedings of the NATO RTO-MP-AVT-147 Computational Uncertainty in Military
  Vehicle design symposium, 2007.

\bibitem{matthies2007}
{\sc H.~G. Matthies}, {\em Uncertainty quantification with stochastic finite
  elements}, in Encyclopedia of Computational Mechanics, E.~Stein, R.~de~Borst,
  and T.~R.~J. Hughes, eds., John Wiley \& Sons, Chichester, 2007.

\bibitem{ong2008}
{\sc Y.~Ong, K.~Lum, and P.~Nair}, {\em Hybrid evolutionary algorithm with
  hermite radial basis function interpolants for computationally expensive
  adjoint solvers}, Computational Optimization and Applications, 39 (2008),
  pp.~97--119.

\bibitem{Press2007}
{\sc W.~H. Press}, {\em Numerical recipes 3rd edition: The art of scientific
  computing}, Cambridge university press, 2007.

\bibitem{radovic1996}
{\sc I.~Radovi\'{c}, I.~M. Sobol, and R.~F. Tichy}, {\em Quasi-monte carlo
  methods for numerical integration: Comparison of different low discrepancy
  sequences}, Monte Carlo Methods and Applications, 2 (1996), pp.~1--14.

\bibitem{ResminiPhD}
{\sc A.~Resmini}, {\em Sensitivity analysis for {RANS} flows. Ph.D. thesis},
  PhD thesis, ONERA-Univ. Paris VI, 2015.

\bibitem{Resmini16}
{\sc A.~Resmini, J.~Peter, and D.~Lucor}, {\em Sparse grids-based stochastic
  approximations with applications to aerodynamics sensitivity analysis},
  International Journal for Numerical Methods in Engineering, 106 (2016),
  pp.~32--57.
\newblock nme.5005.

\bibitem{Schaback1995}
{\sc R.~Schaback}, {\em Error estimates and condition numbers for radial basis
  function interpolation}, Advances in Computational Mathematics, 3 (1995),
  pp.~251--264.

\bibitem{VSCS12}
{\sc C.~Schillings and V.~Schulz}, {\em On the influence of robustness measures
  on shape optimization with stochastic uncertainties}, Optimization and
  Engineering,  (2012), pp.~1--40.

\bibitem{Schillings}
{\sc V.~Schulz and C.~Schillings}, {\em On the nature and treatment of
  uncertainties in aerodynamic design}, AIAA Journal, 47 (2009), pp.~646 --
  654.

\bibitem{snyder2000}
{\sc W.~C. Snyder}, {\em Accuracy estimation for quasi-monte carlo
  simulations}, Mathematics and Computers in Simulation, 54 (2000),
  pp.~131--143.

\bibitem{Wendland2005}
{\sc H.~Wendland}, {\em Scattered data approximation}, Cambridge University
  Press, 2005.

\bibitem{wiener1938}
{\sc N.~Wiener}, {\em The homogeneous chaos}, Amer. J. Math., 60 (1938),
  pp.~897--936.

\bibitem{xiu2010numerical}
{\sc D.~Xiu}, {\em Numerical methods for stochastic computations: a spectral
  method approach}, 2010.

\bibitem{zimmermann2010}
{\sc R.~Zimmermann}, {\em Asymptotic behavior of the likelyhood function of
  covariance matrices of spatial gaussian processes}, J. App. Math.,  (2010).
\newblock Article ID 494070.

\end{thebibliography}

\end{document}